\documentclass[reqno,11pt]{amsart}

\usepackage[a4paper,left=23mm,right=23mm,top=30mm,bottom=30mm,marginpar=25mm]{geometry}
\usepackage{amsmath}
\usepackage{amssymb}
\usepackage{amsthm}
\usepackage{xcolor}
\usepackage{amscd}
\usepackage{mathtools}
\usepackage{url}
\usepackage[ansinew]{inputenc}
\usepackage{color}
\usepackage[final]{graphicx}
\usepackage{enumitem}
\usepackage{esint}
\usepackage{bbm}
\usepackage{subfigure}

\usepackage{tikz}
\usetikzlibrary{arrows,shapes,automata,backgrounds,petri,positioning}
\usetikzlibrary{decorations.pathmorphing}
\usetikzlibrary{decorations.shapes}
\usetikzlibrary{decorations.text}
\usetikzlibrary{decorations.fractals}
\usetikzlibrary{decorations.footprints}
\usetikzlibrary{shadows}
\usetikzlibrary{calc}
\usetikzlibrary{spy}

\usepackage[noadjust]{cite}

\usepackage[colorlinks=true,linktocpage=true, pagebackref=false, citecolor=black,linkcolor=black]{hyperref}

\definecolor{royalblue(web)}{rgb}{0.25, 0.41, 0.88}
\hypersetup{
     colorlinks=true,
     linkcolor=royalblue(web),
     filecolor=magenta,      
     urlcolor=cyan,
     pdfpagemode=FullScreen,
    }

\renewcommand{\epsilon}{\varepsilon}

\numberwithin{equation}{section}

\newtheoremstyle{thmlemcorr}{10pt}{10pt}{\itshape}{}{\bfseries}{.}{10pt}{{\thmname{#1}\thmnumber{ #2}\thmnote{ (#3)}}}
\newtheoremstyle{thmlemcorr*}{10pt}{10pt}{\itshape}{}{\bfseries}{.}\newline{{\thmname{#1}\thmnumber{ #2}\thmnote{ (#3)}}}
\newtheoremstyle{defi}{10pt}{10pt}{\itshape}{}{\bfseries}{.}{10pt}{{\thmname{#1}\thmnumber{ #2}\thmnote{ (#3)}}}
\newtheoremstyle{remexample}{10pt}{10pt}{}{}{\bfseries}{.}{10pt}{{\thmname{#1}\thmnumber{ #2}\thmnote{ (#3)}}}
\newtheoremstyle{ass}{10pt}{10pt}{}{}{\bfseries}{.}{10pt}{{\thmname{#1}\thmnumber{ A#2}\thmnote{ (#3)}}}

\theoremstyle{thmlemcorr}
\newtheorem{theorem}{Theorem}
\numberwithin{theorem}{section}
\newtheorem{lemma}[theorem]{Lemma}
\newtheorem{corollary}[theorem]{Corollary}
\newtheorem{proposition}[theorem]{Proposition}

\newtheorem{conjecture}[theorem]{Conjecture}

\theoremstyle{thmlemcorr*}
\newtheorem{theorem*}{Theorem}
\newtheorem{lemma*}[theorem]{Lemma}
\newtheorem{corollary*}[theorem]{Corollary}
\newtheorem{proposition*}[theorem]{Proposition}
\newtheorem{problem*}[theorem]{Problem}
\newtheorem{conjecture*}[theorem]{Conjecture}

\theoremstyle{defi}
\newtheorem{definition}[theorem]{Definition}

\theoremstyle{remexample}
\newtheorem{remarkx}[theorem]{Remark}
\newenvironment{remark}
  {\pushQED{\qed}\remarkx}
  {\popQED\endremarkx}
\newtheorem{examplex}[theorem]{Example}
\newenvironment{example}
  {\pushQED{\qed}\examplex}
  {\popQED\endexamplex}

\theoremstyle{ass}

\newcommand{\Dcal}{\mathcal{D}}

\newcommand{\Hcal}{\mathcal{H}}
\newcommand{\Ical}{\mathcal{I}}

\newcommand{\Pcal}{\mathcal{P}}

\newcommand{\Scal}{\mathcal{S}}

\DeclareMathOperator{\dist}{dist}

\DeclareMathOperator{\supp}{supp}

\newcommand{\normlr}[1]{\left\|#1\right\|}

\newcommand{\abslr}[1]{\left|#1\right|}

\newcommand{\N}{\mathbb{N}}
\newcommand{\R}{\mathbb{R}}

\newcommand{\C}{\mathbb{C}}

\newcommand{\Z}{\mathbb{Z}}

\DeclareMathOperator{\Div}{div}

\def\XXint#1#2#3{{\setbox0=\hbox{$#1{#2#3}{\int}$}
\vcenter{\hbox{$#2#3$}}\kern-.5\wd0}}

\definecolor{KUblue}{RGB}{0,102,153}

\DeclarePairedDelimiter\abs{\lvert}{\rvert}
\DeclarePairedDelimiter{\norm}{\lVert}{\rVert}

\renewcommand{\phi}{\varphi}

\def\XXint#1#2#3{{\setbox0=\hbox{$#1{#2#3}{\int}$}
     \vcenter{\hbox{$#2#3$}}\kern-.5\wd0}}

\makeatletter
\g@addto@macro\bfseries{\boldmath}
\makeatother

\newcommand{\mres}{\mathbin{\vrule height 1.4ex depth 0pt width
0.13ex\vrule height 0.13ex depth 0pt width 1.0ex}}

\newcommand\restr[2]{{
  \left.\kern-\nulldelimiterspace
  #1 
  \vphantom{\big|} 
  \right|_{#2}
  }}

\usepackage{esint}

\usepackage{appendix}

\DeclareMathOperator{\Tr}{Tr}
\newcommand{\Trt}{\widetilde{\mathrm{Tr}}\,}

\newcommand{\ds}{(-\Delta)^{s/2}}

\renewcommand{\O}{\Omega}

\title[Intrinsic characterizations of Bessel potential spaces on domains]{Intrinsic characterizations of Bessel potential spaces on domains using fractional gradients}

\author{Jos\'e C. Bellido}
\address{ETSI Industrial, Departamento de Matem\'aticas, Universidad de Castilla-La Mancha, Campus Universitario s/n, 13071-Ciudad Real, Spain}
\email{josecarlos.bellido@uclm.es}

\author{Guillermo Garc\'ia-S\'aez}
\address{ETSI Industrial, Departamento de Matem\'aticas, Universidad de Castilla-La Mancha, Campus Universitario s/n, 13071-Ciudad Real, Spain}
\email{guillermo.garciasaez@uclm.es}

\author{Hidde Sch\"{o}nberger}
\address{Department of Applied Mathematics, University of Twente, Drienerlolaan 5, 7522 NB Enschede, The Netherlands}
\email{hidde.schonberger@utwente.nl}

\begin{document}
	
\maketitle

\thispagestyle{empty}
\begin{abstract}  
We establish new characterizations of Bessel potential spaces on domains using fractional gradients that depend only on values within the domain, specifically, the regional fractional gradient and two extended fractional gradients. Our main result characterizes membership in the Bessel potential space by requiring that the function and its fractional gradient belong to the corresponding Lebesgue space on the domain, together with suitable control of the growth near the boundary. This growth is described by a new type of fractional Hardy inequality that also applies to functions with non-zero trace. In certain parameter regimes, the characterizations based on the extended fractional gradients require no separate control of the growth near the boundary. Moreover, the extended fractional gradients provide equivalent norms on the Bessel potential spaces for almost all parameter values. We use this to prove fractional Poincar\'e-Wirtinger inequalities, which provide a basis for studying new classes of fractional partial differential equations on domains.

\vspace{8pt}


\vspace{8pt}


\end{abstract}
	
	\section{Introduction}

Bessel potential spaces were introduced in the early 1960s as a
framework for the study of differentiability of noninteger order related to operators arising in the modern theory of
partial differential equations. The works of Aronszajn and Smith
\cite{ArSmith61}, Calder\'on \cite{Cal61}, and Aronszajn, Mulla and
Szeptycki \cite{ArMuSz63} established the foundations of this theory
and its relation to classical Sobolev spaces. The passage from the
whole space to domains was already an important question in these
developments. In particular, Adams, Aronszajn and Smith
\cite{AdArSm67} studied restrictions of Bessel potentials to domains,
intrinsic descriptions of the resulting spaces, and their connection
with simultaneous extension operators.

Interpolation theory provided another fundamental perspective on
fractional differentiability. The work of Lions and Magenes
\cite{LiMa61} connected Hilbert interpolation spaces with elliptic
boundary value problems, while the complex method developed by
Calder\'on \cite{Cal64} supplied a general framework for constructing
intermediate spaces. For $s \in (0,1)$ and $p \in (1,\infty)$, Bessel potential
spaces satisfy
\[
H^{s,p}(\mathbb R^n)
=
[L^p(\mathbb R^n),W^{1,p}(\mathbb R^n)]_s
\]
with equivalent norms. 

The definition of Bessel potential spaces on a domain is usually
given by restriction from the whole space as
\[
H^{s,p}(\Omega)
=
\left\{
\bar{u}|_\Omega:\bar{u}\in H^{s,p}(\mathbb R^n)
\right\} \quad \text{with} \quad
\|u\|_{H^{s,p}(\Omega)}
=
\inf_{\bar{u}|_\Omega=u}
\|\bar{u}\|_{H^{s,p}(\mathbb R^n)}.
\] 
On bounded Lipschitz domains, restriction and
extension operators give the corresponding identification
\[
H^{s,p}(\Omega)
=
[L^p(\Omega),W^{1,p}(\Omega)]_s;
\]
see \cite{BGS_Compactness,Cal64,Ryc99} for the interpolation and extension frameworks.
A recent treatment of the whole-space identification and its
consequences for continuous embeddings is given in
\cite{BGS_Complex}. The distinction between real and complex
interpolation is relevant here: Bessel potential spaces and
Sobolev--Slobodeckij spaces provide different scales of fractional
regularity, which coincide in the Hilbert case but differ in general
when $p\ne2$.

Although the definition of Bessel potential spaces on domains is well suited to extension and
interpolation arguments, a description using only information within
$\Omega$ requires additional analysis near the boundary. Within the
Besov and Triebel--Lizorkin scales, Rychkov obtained intrinsic
characterizations by local means on smooth domains \cite{Ryc98},
and subsequently constructed a universal extension operator and
intrinsic descriptions on Lipschitz domains \cite{Ryc99}. More recent
developments include characterizations of generalized Besov,
Triebel--Lizorkin and Morrey-type spaces on special Lipschitz domains
\cite{Yao23}, descriptions by averaged differences on uniform domains
and their application to composition and inverse regularity
\cite{Prats24}, and derivative norm equivalences together with new
estimates for Rychkov's extension operator \cite{ShY24}. These
results show how the geometry of the domain enters the
characterization of fractional smoothness. Despite these results, a simple characterization as for the Sobolev-Slobodeckij spaces on domains \cite{DPV12}
\[
W^{s,p}(\O)=\left\{u \in L^p(\O)\,:\, \left(\int_{\O}\int_{\O} \frac{|u(y)-u(x)|^p}{|y-x|^{n+sp}}\,dy\,dx\right)^{1/p} < \infty \right\},
\]
which opens the way to study fractional differential equations on domains, is still absent for $H^{s,p}(\O)$.

To fill this gap, we start from the recent differential interpretation of Bessel potential spaces arising
from the Riesz fractional gradient studied by Shieh and Spector in
\cite{ShS2015}. For a smooth compactly supported function, this
operator is given by
\[
D^s u(x)
=
c_{n,s}\int_{\mathbb R^n}
\frac{u(y)-u(x)}{|y-x|^{n+s}}
\frac{y-x}{|y-x|}\,dy,
\qquad s \in (0,1),
    \]
where $c_{n,s}$ is a normalization constant. In
\cite[Theorem~1.7]{ShS2015}, the completion of
$C_c^\infty(\mathbb R^n)$ with respect to
\[
\|u\|_{L^p(\mathbb R^n)}
+
\|D^s u\|_{L^p(\mathbb R^n;\mathbb R^n)}
\]
is identified with $H^{s,p}(\mathbb R^n)$, with equivalent norms, while 
the corresponding distributional formulation is developed further
in \cite{BrCaCoSt22}. This connects potential theory with a
fractional differential structure analogous to that of classical
Sobolev spaces, and it is well suited for problems in the calculus of variations and partial differential equations, see~e.g.,~\cite{ABS25, BeCuMC, Borthagaray2025, Camilo2026, camposI,CamposGarciaSaez2026, Carrero2025, Gambera2024, GarciaSaez2026, KrS22, Lo2023, Schikorra2018, Silhavy2024, ShS2018}, among many others. For functions defined only on a bounded domain, however, the Riesz
fractional gradient cannot immediately be applied. Choosing
an extension of the function makes the operator available, but its restriction to
$\Omega$ generally depends on that choice. Extension by zero is
particularly useful for supported Bessel spaces, but it does not
preserve $H^{s,p}$ regularity for all parameters. In the presence of
traces, its regularity entails a zero boundary restriction and therefore
does not cover all of $H^{s,p}(\Omega)$. Thus, the whole-space
differential characterization does not directly produce an operator determined by the function
on $\Omega$ alone.

\subsection*{Contributions of the paper} In this paper, we investigate fractional differential descriptions
of Bessel potential spaces that use only information on the given
domain. Our starting point is the regional fractional gradient
\[
D^s_\Omega u(x)
=
c_{n,s}\int_\Omega
\frac{u(y)-u(x)}{|y-x|^{n+s}}
\frac{y-x}{|y-x|}\,dy,
\qquad x\in\Omega,
\]
initially defined for smooth functions and then interpreted weakly. Its natural domain is the space
\[
\Hcal^{s,p}(\O)
:=\{u\in L^p(\O):D^s_\O u\in L^p(\O;\R^n)\},
\]
endowed with the graph norm
$\|u\|_{L^p(\O)}+\|D^s_\O u\|_{L^p(\O;\R^n)}$.
The regional construction requires neither values in an exterior
collar nor a horizon parameter, which distinguishes it from the
finite-horizon gradients introduced in \cite{BeCuMC23}; for the closely related regional fractional Laplacian, we refer to \cite{LiP20} and the references therein. The central question motivating our analysis is the following.

\begin{conjecture}\label{conj:equality Bessel and regional}
Let $\O\subset\R^n$ be a smooth and bounded domain,
$s \in (0,1)$ and $p \in (1,\infty)$. Then,
\[
\Hcal^{s,p}(\O)=H^{s,p}(\O)
\]
with equivalent norms.
\end{conjecture}

While this is still an open problem in its full generality, our results equate the interior regularity of both spaces and identify the boundary growth behavior that recovers the global Bessel regularity. Precisely, for every bounded open set $\O \subset \R^n$ and
the full range $s \in (0,1)$, $p \in (1,\infty)$, we prove
\begin{equation}\label{eq:intembeddings}
H^{s,p}(\O)\hookrightarrow\Hcal^{s,p}(\O)
\hookrightarrow L^p(\O)\cap H^{s,p}_{\mathrm{loc}}(\O);
\end{equation}
see Theorems~\ref{th:emb1} and~\ref{th:regional-local-bessel}.
The first embedding follows from fractional commutator estimates
of Li \cite{li2019}. The second is obtained by localizing the
Riesz gradient of the zero-extension. 
In parallel, we prove that a vanishing regional fractional gradient implies constancy; it holds either for functions in $C^2(\overline{\O})$ or functions in $L^1(\O)$ when $\O$ is convex. This is in stark contrast with the truncated fractional gradients, where the values in the collar enable non-constant functions with zero truncated fractional gradient \cite{KrS24}. These results address a necessary step towards a regional Poincar\'e-Wirtinger inequality. 

The remaining issue for identifying the spaces concerns control up to the boundary, which turns out to be closely connected to the following two boundary estimates that hold on the Bessel potential space: The fractional Hardy inequality
\[
\|u/d^s\|_{L^p(\O)} \lesssim \|u\|_{H^{s,p}(\R^n)} \quad \text{for all $u \in H^{s,p}_0(\O)$,}
\]
where $d(x):=\dist(x,\O^c)$, and, when $sp>1$, the estimate (Corollary~\ref{cor:hardy})
\[
\left\|
\int_{\partial\O}\nu_\O(y)
\frac{\Tr u(y)-u(\cdot)}{|y-\cdot|^{n+s-1}}
\,d\Hcal^{n-1}(y)
\right\|_{L^p(\O;\R^n)}
\lesssim\|u\|_{H^{s,p}(\O)} \quad \text{for all $u \in H^{s,p}(\O)$,}
\]
where $\nu_\O$ is the inner unit normal; the latter estimate incorporates
non-zero boundary values into a fractional Hardy-type condition, and implies the fractional Hardy inequality when $\Tr u =0$. This inequality seems to be new in the literature, and we prove it using the boundedness of $D^s_\O$ on $H^{s,p}(\O)$ from \eqref{eq:intembeddings}. Our main result (Theorem~\ref{th:char}) shows that for $\O \subset \R^n$ smooth and bounded and $p\ne 2/(1+s)$, these Hardy inequalities provide the missing piece for the identification. Precisely, it holds that $u \in H^{s,p}_0(\O)$ if and only if
\[
u \in \Hcal^{s,p}(\O) \ \ \text{and} \ \ u/d^s \in L^p(\O),
\]
while, for $sp>1$, it holds that $u \in H^{s,p}(\O)$ if and only if
\[
u \in \Hcal^{s,p}(\O) \quad \text{and} \quad \int_{\partial \Omega}\nu_\Omega(y) \frac{g(y)-u(\cdot)}{|y-\cdot|^{n+s-1}}\,d\Hcal^{n-1}(y) \in L^p(\O;\R^n),
\]
for some $g \in W^{s-1/p,p}(\partial \O)$; the boundary datum $g$ is then necessarily the trace of $u$. Since $H^{s,p}_0(\O) \cong H^{s,p}(\O)$ when $sp<1$, the first characterization is also applicable to $H^{s,p}(\O)$ in that range. The cause of the exceptional value $p=2/(1+s)$ is explained in Remark~\ref{rem:oneinclusion}\,c). 

To obtain these characterizations, we introduce two extended
fractional gradients, which also only depend on values within $\O$:
\[
\bar G^s_\O u=D^s(\mathbbm{1}_\O u)|_\O
\quad \text{and} \quad
G^s_\O u=I_{1-s}(\mathbbm{1}_\O\nabla u)|_\O,
\]
with $I_{1-s}$ the Riesz potential of order $1-s$. The difference between them is a Riesz potential of the boundary values,
as expressed in \eqref{eq:identityG}; this shows, in particular, that the latter operator is only suitable in the range $sp>1$.
A careful application of the theory on transmission spaces developed by Grubb \cite{Gru15} allows us to describe the weak domains
\[
\{u \in L^p(\O) \,:\, \bar{G}^s_\O u \in L^p(\O;\R^n)\} \quad \text{and} \quad \{u \in L^p(\O) \,:\, G^s_\O u \in L^p(\O;\R^n)\}
\]
and the parameter ranges in which they coincide with
the corresponding Bessel potential spaces; see Theorem~\ref{th:tildechar}. This provides another characterization of the Bessel potential spaces in terms of fractional gradients that complements Theorem~\ref{th:char}.
Surprisingly, even for different parameter ranges, the norms induced by the extended fractional gradients provide equivalent norms on the Bessel potential spaces (Proposition~\ref{prop:equivnorm}). In particular, for smooth bounded
domains and $sp>1$, the norm of $H^{s,p}(\O)$ is equivalent to
\[
\|u\|_{L^p(\O)}+\|G^s_\O u\|_{L^p(\O;\R^n)}.
\]
This leads to Poincar\'e-Wirtinger
inequalities for $\bar{G}^s_\O$ and $G^s_\O$ in Corollary~\ref{cor:poincwirt}, which open up the possibility of studying fractional partial differential equations and variational problems.
The description of the full weak domains is a separate issue
from these norm equivalences: Boundary singularities are not ruled out by the extended fractional gradients, which results in functions that cannot be approximated by smooth functions, cf.~Remark~\ref{rem:singularities}. This is closely related to the recent results for truncated versions of $\bar{G}^s_\O$ in the case $p=2$, see~\cite{FLP26}.

\subsection*{Outline} The paper is organized as follows: Section~\ref{sec:preliminaries} contains preliminaries on Bessel potential spaces and the transmission spaces from \cite{Gru15}. In Section~\ref{sec:regional}, we introduce the regional fractional spaces $\Hcal^{s,p}(\O)$ with Section~\ref{sec:nullspace} containing the description of the null-space of $D^s_\O$ and Section~\ref{sec:connection} the connections between $H^{s,p}(\O)$ and $\Hcal^{s,p}(\O)$. Finally, Section~\ref{sec:extended} is devoted to the extended fractional gradients $\bar{G}^s_\O$ and $G^s_\O$ and their associated spaces. This includes the relation with the Bessel potential spaces and the proofs of the equivalence of norms and fractional Poincar\'e-Wirtinger inequalities.

\section{Preliminaries}\label{sec:preliminaries}
\subsection{Notation}
We fix $n\in \mathbb{N}$ for the dimension of our ambient space $\R^n$ and we will denote by $\Omega\subset \R^n$ an open and bounded subset, with additional assumptions mentioned in the respective results. We write $|x|$ for the Euclidean norm of a vector $x=(x_1,\ldots,x_n)\in \R^n$. 
The ball centered at $x\in \R^n$ and with radius $r>0$ is denoted by $B_r(x)=\{y\in \R^n: |x-y|<r\}$. The complement of a set $E\subset \R^n$ is denoted by $E^c=\R^n\setminus{E}$, its closure by $\overline{E}$ and its boundary by $\partial E$. The distance between a point $x\in \R^n$ and a set $E\subset\R^n$ is denoted by $d(x,E)$.

The notation for Sobolev $W^{1,p}$ and Lebesgue $L^p$ spaces is the standard one, as is that of smooth
functions of compact support $C_c^\infty$. We will
indicate the domain of the functions, as in $L^p
(\Omega)$; the target is indicated only if it is not $\R$.  For $u \in L^1(\O)$, we denote its average by
\[
(u)_\O:=\frac{1}{|\O|}\int_{\O} u \,dx.
\]
The notation for the classical gradient in its weak sense is $\nabla$. For a Banach space $X$, we denote its topological dual as $X'$, and the usual {\it strong} convergence in normed spaces is denoted as $\to$. For a distribution $T \in \Dcal'(\R^n)$, we denote by $\langle T,\varphi\rangle$ the usual pairing for $\varphi\in C_c^\infty(\R^n)$.

Our convention for the Fourier transform of functions $f\in L^1(\R^n)$ is 
\[
\widehat{f}(\xi)=\int_{\R^n}f(x)e^{-2\pi i x\cdot \xi}\,dx \quad \text{for $\xi\in \R^n$}.
\]
This definition is extended by continuity and duality to other function and distribution spaces as usual in the theory of function spaces. The inverse Fourier transform of $f$, denoted $f^\vee$, corresponds to $x \mapsto \widehat{f}(-x)$.
We denote by $\mathcal{S}(\R^n)$ the Schwartz space and $\mathcal{S}'(\R^n)$ the space of tempered distributions. More details can be found in the classical text \cite{Gra14a}. 

\subsection{Bessel potential spaces and fractional gradients}
We recall that the Bessel potential spaces for $s \in \R$ and $p \in (1,\infty)$ are given by
\[
H^{s,p}(\R^n):=\left\{u \in \Scal'(\R^n)\,:\, \left(\langle \cdot \rangle^{s}\widehat{u}\right)^{\vee} \in L^p(\R^n)\right\} \quad \text{with} \quad \norm{u}_{H^{s,p}(\R^n)}:=\norm{\left(\langle \cdot \rangle^{s}\widehat{u}\right)^{\vee}}_{L^p(\R^n)},
\]
where $\langle \xi \rangle :=(1+|\xi|^2)^{1/2}$, sometimes referred to as the \textit{Japanese bracket}, is a smooth function that behaves like $|\xi|$ for large $\xi$. We also define for an open set $\Omega \subset \R^n$ the closed subspaces $H^{s,p}_0(\Omega) \subset H^{s,p}(\R^n)$ as
\[
H^{s,p}_0(\Omega):=\{u \in H^{s,p}(\R^n) \,:\, \supp(u) \subset \overline{\Omega}\},
\]
where $\supp(u)$ is understood in the sense of distributions. When $s \in (0,1)$ and $u \in L^p(\O)$ we sometimes use a slight abuse of notation and write $u \in H^{s,p}_0(\O)$, when we mean $\mathbbm{1}_\O u \in H^{s,p}_0(\O)$. The restricted space is defined as
\[
H^{s,p}(\Omega):=\{u \in \Dcal'(\Omega) \,:\, \exists\bar{u} \in H^{s,p}(\R^n) \ \text{with} \ \bar{u}|_\Omega = u\}
\]
endowed with the quotient norm
\[
\norm{u}_{H^{s,p}(\Omega)}:=\inf_{\bar{u} \in H^{s,p}(\R^n), \ \bar{u}|_\Omega = u} \norm{\bar{u}}_{H^{s,p}(\R^n)}.
\]
We note that $H^{s,p}(\R^n)$ and $H^{-s,p'}(\R^n)$ are dual spaces with $1/p+1/p'=1$, and, when $C_c^{\infty}(\Omega)$ is dense in $H^{s,p}_0(\Omega)$, the same holds for $H^{s,p}_0(\Omega)$ and $H^{-s,p'}(\Omega)$. Such a density result holds, for example, for bounded Lipschitz domains, while $C^{\infty}(\overline{\O})$ is always dense in $H^{s,p}(\O)$ by extension. In addition, for Lipschitz domains $\O$, $s \in (0,1)$ and $sp>1$ there exists a bounded surjective trace operator
\begin{equation}\label{eq:trace}
\Tr:H^{s,p}(\O) \to W^{s-1/p,p}(\partial \O),
\end{equation}
where
\[
W^{s-1/p,p}(\partial \O)=\left\{u \in L^p(\partial \O) \,:\, \int_{\partial \O}\int_{\partial\O} \frac{|u(y)-u(x)|^{p}}{|y-x|^{n-2+sp}}\,d\Hcal^{n-1}(y)\,d\Hcal^{n-1}(x) < \infty \right\},
\]
see~e.g.~\cite[Theorem~1.1]{Gau23}. We also make frequent use of the following estimate
\begin{equation}\label{eq:hspnorm}
    \|u\|_{H^{s,p}(\O)} \lesssim \|u\|_{L^p(\O)}+\|\nabla u\|_{H^{s-1,p}(\O)} \quad \text{for all $u \in H^{s,p}(\O)$,}
\end{equation}
for Lipschitz domains and $s \in (0,1)$, see~\cite{ShY24}.

Bessel potential spaces generalize the classical Sobolev spaces $W^{k,p}$ for integer $k\in \Z$ to every real number $s\in \R$, and can be obtained by complex interpolation methods between $L^p$ and $W^{k,p}$ (see \cite{BGS_Complex}). For $H^{s,p}$-spaces the analogous continuous and compact embeddings hold depending on the value of the parameter $sp-n$, as well as the structural properties that one could expect. In particular, we want to recall the relationship with the extensively used \textit{Sobolev-Slobodeckij} or \textit{Gagliardo spaces} $W^{s,p}$, defined by $$W^{s,p}(\Omega):=\Bigg\{u\in L^p(\Omega): \frac{|u(y)-u(x)|}{|y-x|^{\frac{n}{p}+s}}\in L^p(\Omega\times \Omega)\Bigg\}$$ with $$\norm{u}_{W^{s,p}}:=\left(\norm{u}_p^p+[u]_{W^{s,p}}^p\right)^\frac{1}{p} \quad \text{where} \quad [u]_{W^{s,p}}:=\left(\int_{\Omega}\int_{\Omega}\frac{|u(y)-u(x)|^p}{|y-x|^{n+sp}}\,dx\,dy\right)^{1/p},$$ is the Gagliardo seminorm. Those spaces are of real interpolation nature, and exploiting the relationship between both interpolation methods in \cite{BGS_Complex} the following result was proved. 
\begin{theorem}[Bessel versus Gagliardo]\label{Contiguity}
    Let $\Omega\subset \R^n$ be a Lipschitz domain, $0<s_0<s<s_1<1$ and $p \in (1,\infty)$. Then, \begin{itemize}
        \item[(i)] $W^{s,2}(\Omega)=H^{s,2}(\Omega)$ with equivalence of the norms. \smallskip
        \item[(ii)] $H^{s_1,p}(\Omega)\xhookrightarrow{}W^{s,p}(\Omega)\xhookrightarrow{}H^{s_0,p}(\Omega).$ \smallskip
        \item[(iii)] It holds that $$\begin{cases}
        W^{s,p}(\Omega)\xhookrightarrow{}H^{s,p}(\Omega),&1<p\leq 2,\\
        H^{s,p}(\Omega)\xhookrightarrow{}W^{s,p}(\Omega),&2\leq p<\infty,
    \end{cases}$$ with strict inclusions unless $p=2$.
    \end{itemize}
\end{theorem}
In \cite{ShS2015,ShS2018}, the 
\textit{Riesz fractional gradient and divergence} of $u\in C_c^\infty(\R^n)$ and $v\in C_c^\infty(\R^n;\R^n)$, respectively, were introduced as 
\begin{align*}
D^su(x)&:=c_{n,s}\int_{\R^n}\frac{u(y)-u(x)}{|y-x|^{n+s}}\frac{y-x}{|y-x|}\,dy \quad \text{for $x\in \R^n$},\\
\operatorname{div}^sv(x)&:=c_{n,s}\int_{\R^n}\frac{v(y)-v(x)}{|y-x|^{n+s}}\cdot\frac{y-x}{|y-x|}\,dy \quad \text{for $x\in \R^n$},
\end{align*}
where $c_{n,s}>0$ is some normalization constant.

The following proposition recalls the main properties of this operator that we will require for our development, we refer to \cite{BeCuMC21,KrS22,ShS2015,ShS2018} for detailed proofs.
\begin{proposition}\label{DsProps}
     Let $s\in (0,1)$, $p \in (1,\infty)$, $u\in C_c^\infty(\R^n)$ and $v\in C_c^\infty(\R^n;\R^n)$. Then, the following hold:
     \begin{itemize}
         \item[(i)] The space $\{u\in L^p(\R^n): D^su\in L^p(\R^n;\R^n)\}$ endowed with the norm $$\norm{\cdot}_{s,p}:=\norm{\cdot}_p+\norm{D^s\cdot}_p$$ 
         and $\overline{C^\infty_c(\R^n)}^{\norm{\cdot}_{s,p}}$ both coincide with $H^{s,p}(\R^n)$ with equivalence of norms. \smallskip
         \item[(ii)] We have the identities $$D^su=I_{1-s}\left(\nabla u\right)=\nabla (I_{1-s}u) \quad \text{and} \quad \operatorname{div}^sv=I_{1-s}(\operatorname{div}v)=\operatorname{div}(I_{1-s}v),$$ where $I_{s}u$ denotes the \textit{Riesz potential} of $u$, defined as a convolution operator with kernel $\gamma_{n,s}^{-1}|\cdot|^{-(n-s)}$, for some normalization constant $\gamma_{n,s}>0$. \smallskip
         \item[(iii)] The following integration by parts formula holds:$$\int_{\R^n}D^su(x)\cdot v(x)\,dx=-\int_{\R^n}u(x)\operatorname{div}^sv(x)\,dx.$$ \smallskip
         \item[(iv)] The Fourier transforms of $D^su$  and  $\operatorname{div}^sv$  are $$\widehat{D^su}(\xi)=\frac{2\pi i\xi}{|2\pi\xi|^{1-s}}\widehat{u}(\xi) \quad \text{and} \quad \widehat{\operatorname{div}^sv}(\xi)=\frac{2\pi i\xi}{|2\pi\xi|^{1-s}}\cdot \widehat{v}(\xi) \quad \text{for $\xi \in \R^n \setminus \{0\}$}.$$ \smallskip
         \item[(v)] We have the identity $-\operatorname{div}^sD^su=(-\Delta)^su$ and $D^su=\mathcal{R}(-\Delta)^{s/2}u$, where $\mathcal{R}$ is the vector-valued Riesz transform.
     \end{itemize}
\end{proposition}

\subsection{Transmission spaces}
We recall here the transmission spaces introduced in \cite{Gru15} for studying the regularity of solutions to pseudo-differential equations.

To define the transmission spaces, we assume $\Omega \subset \R^n$ is a smooth and bounded domain, and use for $\mu \in \R$ the order reducing operators $\Lambda^{(\mu)}_{+}$ from \cite[Theorem~1.3]{Gru15}, which are homeomorphisms
\begin{equation}\label{eq:orderreducing}
\Lambda^{(\mu)}_{+}: H^{s,p}_0(\Omega) \to H^{s-\mu,p}_0(\Omega) \quad \text{for all $s \in \R$,}
\end{equation}
among other properties. The transmission spaces are now defined as follows.
\begin{definition}[Transmission spaces]\label{def:transmission} Let $\mu,s \in \R$, $p \in (1,\infty)$ with $s > \mu -1/p'$. Then, the transmission space $H^{\mu(s)}_p(\overline{\Omega})$ consists of all distributions $u \in \Scal'(\R^n)$ with $\supp (u) \subset \overline{\Omega}$ such that $\Lambda^{(\mu)}_{+}u \in \bigcup_{\epsilon >0} H^{\epsilon-1/p',p}_0(\Omega)$ and
\[
\|u\|_{H^{\mu(s)}_p(\overline{\Omega})}:=\|(\Lambda^{(\mu)}_{+}u)|_{\Omega}\|_{H^{s-\mu,p}(\Omega)} < \infty.
\]
\end{definition}
Note that since $\supp(\Lambda^{(\mu)}_{+}u) \subset \overline{\Omega}$, we find that for $s <\mu +1/p$ actually $\Lambda^{(\mu)}_{+}u \in H^{s-\mu,p}_0(\Omega)$, since extension by zero is bounded in that case. Hence, by \eqref{eq:orderreducing}, we infer that $H^{\mu(s)}_p(\overline{\Omega})=H^{s,p}_0(\Omega)$ for $s<\mu+1/p$. On the other hand, if $\mu+1/p <s<\mu + 1/p+1$ (and $\mu>-1$), then we find by \cite[Theorem~5.4]{Gru15} that
\[
H^{\mu(s)}_p(\overline{\Omega}) \subset d^\mu H^{s-\mu,p}(\Omega) + H^{s,p}_0(\Omega),
\]
with $d = \dist(\cdot,\Omega^c)$. Moreover, $u/d^\mu$ has a well-defined trace at $\partial \Omega$; if this trace vanishes then $u \in H^{s,p}_0(\Omega)$. This essentially shows that $H^{\mu(s)}_p(\overline{\Omega})$ is a space with smoothness governed by the parameter $s$, and boundary behavior regulated by $\mu$.

Next, we will present the main regularity result involving the transmission spaces. We refer to \cite{Gru15} for the precise definitions of some of the properties of pseudo-differential operators; we will only apply it to the Riesz potential.

\begin{theorem}[{\cite[Theorem~4.4]{Gru15}}]\label{th:regularity}
Let $\Pcal$ be a classical elliptic pseudo-differential operator of order $m \in \R$ and of type $\mu \in \R$ relative to $\Omega$ and with factorization index $\mu$. Let $s>\mu - 1/p'$ and $u \in \bigcup_{\epsilon >0} H^{\epsilon+\mu-1/p',p}_0(\Omega)$. Then, it holds that
\[
\Pcal u|_{\Omega} \in H^{s-m,p}(\Omega) \quad \Rightarrow \quad u \in H^{\mu(s)}_p(\overline{\Omega}).
\]  
\end{theorem}
We want to apply this result to the Riesz potential $\Pcal = I_{1-s}$ with $s \in (0,1)$. The symbol of the Riesz potential is given by
\[
p(\xi):=\abs{2\pi \xi}^{s-1}.
\]
This clearly shows that $I_{1-s}$ has order $s-1$, is elliptic, and, since the symbol is homogeneous, it is classical with decomposition $p_0=p$ and $p_j=0$ for $j \not = 0$, cf.~\cite[page 5]{Gru15}. Additionally, since $p(\xi)=p(-\xi)$, it is clear that $I_{1-s}$ is of type $\mu=\frac{s-1}{2}$ relative to any domain, see~\cite[Definition~2.5]{Gru15}. Finally, because $p$ is real-valued, its factorization index is equal to $\mu$ as well, cf.~\cite[Eq.~(3.4)]{Gru15}. We obtain the following corollary.
\begin{corollary}
    For $s \in (0,1)$ the operator $I_{1-s}$ is a classical elliptic pseudo-differential operator of order $s-1$ and with type and factorization index $\frac{s-1}{2}$.
\end{corollary}
We can use it to prove the following, by combining it with an existence argument as in \cite{Gru22}.
\begin{lemma}\label{le:Rieszisomorphism0}
    Let $s \in (0,1)$ and $q \in (1,\infty)$. Then, it holds that
    \[
    I_{1-s}^{\Omega}:u \mapsto (I_{1-s}u)|_{\Omega} \quad \text{is an isomorphism between $H^{\frac{s-1}{2}(0)}_q(\overline{\Omega})$ and $H^{1-s,q}(\Omega)$}.
    \]
\end{lemma}
\begin{proof}
    The operator is continuous by \cite[Theorem~4.2]{Gru15}. It now suffices to prove bijectivity, by the Banach isomorphism theorem. First, we prove it in the case $q=2$. To this aim, we introduce the bilinear form
    \begin{equation}\label{eq:rieszidentity}
    B[u,v]:=\int_{\R^n} I_{1-s}u \,v \,dx=\int_{\R^n} |2\pi \xi|^{s-1}\widehat{u}(\xi) \overline{\widehat{v}(\xi)}\,d\xi \quad \text{for all $u,v \in C_c^{\infty}(\Omega)$,}
    \end{equation}
    and extend it by density to $B:H^{\frac{s-1}{2}}_0(\Omega) \times H^{\frac{s-1}{2}}_0(\Omega) \to \R$. It is not difficult to check that this bilinear form satisfies the assumptions of the Lax-Milgram lemma. Hence, for each $w \in H^{\frac{s-1}{2}}_0(\Omega)'\cong H^{\frac{1-s}{2}}(\Omega)$, there exists a unique $u \in H^{\frac{s-1}{2}}_0(\Omega)$ with
    \[
    B[u,v] = \langle w,v \rangle \quad \text{for all $v \in H^{\frac{s-1}{2}}_0(\Omega)$,}
    \]
    or, equivalently, $I_{1-s}^{\Omega}u=w$. In particular, this holds if $w \in H^{1-s}(\Omega)$, and by Theorem~\ref{th:regularity}, we conclude that then $$u \in H^{\frac{s-1}{2}(0)}_2(\overline{\Omega}),$$ which proves the desired bijectivity when $q=2$; note that the injectivity also uses the embedding
    \[
    H^{\frac{s-1}{2}(0)}_2(\overline{\Omega}) \hookrightarrow H^{\frac{s-1}{2}(\frac{s-1}{2})}_2(\overline{\Omega})=H^{\frac{s-1}{2}}_0(\O).
    \]
    
    To extend the injectivity to other values of $q$, suppose that $I_{1-s}^{\Omega}u = 0$, then we find by Theorem~\ref{th:regularity} and \cite[Eq.~(4.1)]{Gru15} that 
    \[
    u \in \bigcap_{m \in \R} H^{\frac{s-1}{2}(m)}_{q}(\overline{\Omega})=d^{\frac{s-1}{2}}C^\infty(\overline{\Omega}).
    \]
    In particular, this space does not depend on $q$, and thus we find that $u=0$ by the injectivity in the case $q=2$. For the surjectivity, first consider $q>2$. We note that for any $$w \in H^{1-s,q}(\Omega) \subset H^{1-s}(\Omega),$$ there is a $u \in H^{\frac{s-1}{2}(0)}_2(\overline{\Omega})$ with $I^{\Omega}_{1-s}u = w$. Then, by bootstrapping Theorem~\ref{th:regularity} as in \cite[Theorem~3.4]{Gru22}, we can increase the integrability exponent and conclude $$u \in H^{\frac{s-1}{2}(0)}_q(\overline{\Omega})$$ as desired. For $q<2$, it suffices to prove that the range of $I^{\Omega}_{1-s}$ is closed in $H^{1-s,q}(\Omega)$ since the range contains the dense subset $H^{1-s}(\Omega)$. This property of the range follows by adapting the argument in~\cite[Proposition~3.3]{Gru22} to negative orders, whose changes we now detail. Precisely, using the boundedness 
    \[
    I_{1-s}^{\O}: H^{\frac{s-1}{2}(0)}_q(\overline{\Omega}) \to H^{1-s,q}(\O)
    \]
    from \cite[Theorem~4.2]{Gru15} and the regularity from Theorem~\ref{th:regularity}, we find in the same way using the Banach isomorphism theorem that 
    \[
    \|u\|_{H^{\frac{s-1}{2}(0)}_q(\overline{\Omega})} \lesssim \|I_{1-s}^\O u\|_{H^{1-s,q}(\O)}+\|u\|_{H^{\frac{s-1}{2},q}_0(\O)}.
    \]
    In addition, using the embeddings
    \[
    H^{\frac{s-1}{2}(0)}_q(\overline{\Omega}) \hookrightarrow H^{\frac{s-1}{2}(\frac{s-1}{2}+\varepsilon)}_q(\overline{\Omega})=H^{\frac{s-1}{2}+\varepsilon,q}_0(\O) \hookrightarrow H^{\frac{s-1}{2},q}_0(\O),
    \]
    for some $0<\epsilon <\min\{1/q,(1-s)/2\}$, of which the second embedding is compact, we may conclude that the range of $I_{1-s}^\O$ is closed via an argument by contradiction as in~~\cite[Proposition~3.3]{Gru22}.
\end{proof}
A consequence is the following result, where we introduce the exponents
\[
q_s:=\frac{2}{1+s} \quad \text{and} \quad q_s'=\frac{2}{1-s}.
\]
\begin{lemma}\label{le:Rieszisomorphism}
    Let $s \in (0,1)$ and $p \in (1,\infty)$. If $p>q_s$, then it holds that
    \[
    \text{$I_{1-s}^{\Omega}$ is an isomorphism between $L^{p'}_0(\Omega)$ and $H^{1-s,p'}(\Omega)$}.
    \]
    If also $p>q_s'$, then 
    \[
    \text{$I_{1-s}^{\O}$ is an isomorphism between $H^{s,p'}_0(\Omega)$ and $W^{1,p'}(\Omega)$}.
    \]
\end{lemma}
\begin{proof}
The first statement follows immediately from Lemma~\ref{le:Rieszisomorphism0} together with the fact that $$H^{\frac{s-1}{2}(0)}_{p'}(\overline{\Omega})=L^{p'}_0(\Omega)$$ whenever $p'<q_s' \Leftrightarrow p>q_s$; here, $L^{p'}_0(\O)$ denotes the functions in $L^{p'}(\R^n)$ supported on $\overline{\O}$. For the second, we note that $$H^{\frac{s-1}{2}(s)}_{p'}(\overline{\Omega})=H^{s,p'}_0(\O)$$ when $s+\frac{1-s}{2} < 1/p' \Leftrightarrow p>q_s'$.
\end{proof}

\section{The regional gradient and the space $\Hcal^{s,p}(\O)$}\label{sec:regional}

We introduce the \textit{regional fractional gradient}
\[
D^s_\O u(x):=c_{n,s}\int_{\Omega} \frac{u(y)-u(x)}{\abs{y-x}^{n+s}}\frac{y-x}{\abs{y-x}}\,dy \quad \text{for $x \in \O$},
\]
where $s \in (0,1)$, $\Omega \subset \R^n$ is open and bounded and $u \in C^1(\overline{\O})$. As well, we define the \textit{regional fractional divergence} as
\[
\Div^s_\O v(x):=c_{n,s}\, \lim_{r \downarrow 0} \int_{\Omega \setminus B_r(x)} \frac{v(y)+v(x)}{\abs{y-x}^{n+s}}\cdot \frac{y-x}{\abs{y-x}}\,dy \quad \text{for $x \in \O$},
\] 
for $v \in C_c^{1}(\O;\R^n)$. For these two operators, one has the following integration by parts formula (see~\cite[Theorem 1.4]{MeS15}).
\begin{lemma}[Regional integration by parts]
For $u \in C^1(\overline{\O})$ it holds that $D^s_\O u \in L^1(\O;\R^n)$ and
\[
\int_{\O} D^s_\O u \cdot v \,dx = -\int_{\O} u \Div^s_\O v \,dx \quad \text{for all $v \in C_c^{1}(\O;\R^n)$.}
\]
\end{lemma}
In view of this identity, it is natural to define the regional fractional gradient in a weak sense as follows.
\begin{definition}[Weak fractional gradient]
Let $u \in L^1(\O)$, then we say that $V \in L^1_{\rm loc}(\O;\R^n)$ is its weak regional fractional gradient, written as $V=D^s_\Omega u$, if
\[
\int_{\O} V \cdot v \,dx = -\int_{\O} u \Div^s_\O v \,dx \quad \text{for all $v \in C_c^{\infty}(\O;\R^n)$.}
\]
\end{definition}
Using this, one can naturally define a class of fractional Sobolev spaces.

\begin{definition}[Fractional Sobolev space]
Let $p \in [1,\infty]$, then we define the fractional Sobolev space
\[
\Hcal^{s,p}(\O):=\{u \in L^p(\O) \,:\, D^s_\O u \in L^p(\O;\R^n)\},
\]
endowed with the norm
\[
\norm{u}_{\Hcal^{s,p}(\O)}:=\norm{u}_{L^p(\O)} + \norm{D^s_\O u}_{L^p(\O;\R^n)}.
\]
\end{definition}
We begin by observing that this new class consists of spaces that are complete under their respective norms.
\begin{proposition}\label{prop:banach}
    $\Hcal^{s,p}(\O)$ is a Banach space.
\end{proposition}
\begin{proof}
    Let $(u_j)_j \subset \Hcal^{s,p}(\O)$ be a Cauchy sequence, then we find by completeness of $L^p$ that $u_j \to u$ in $L^p(\Omega)$ and $D^s_\Omega u_j \to U$ in $L^p(\O;\R^n)$. It remains to show that $U=D^s_\Omega u$, which follows from the usual distributional argument. Indeed, for any $\phi \in C_c^{\infty}(\O;\R^n)$, we can compute that
    \begin{align*}
        \int_\O U \cdot \phi \,dx &= \lim_{j \to \infty} \int_\O D^s_\O u_j \cdot \phi \,dx\\
        &=-\lim_{j \to \infty} \int_\O u_j \Div^s_\O \phi \,dx \\
        &=-\int_\O u \Div^s_\O \phi\,dx,
    \end{align*}
    since $\Div^s_\O \phi \in C^\infty(\overline{\O})$. This finishes the proof.
\end{proof}
To highlight the connection with the Riesz fractional gradient, we note that for $u \in C^{\infty}_c(\O)$ 
\begin{equation}\label{eq:connDsobarGso}
D^s_\O u(x) = D^s (\mathbbm{1}_\O u)(x) - u(x)D^s \mathbbm{1}_{\Omega}(x) \quad \text{for $x \in \O$.}
\end{equation}
Indeed, 
\begin{align*}
    D^s(\mathbbm{1}_\O u)(x)&=c_{n,s}\int_{\R^n}\frac{\mathbbm{1}_\O(y)u(y)-\mathbbm{1}_\O(x)u(x)}{|y-x|^{n+s}}\frac{y-x}{|y-x|}\,dy \\
    &=c_{n,s}\int_{\O}\frac{u(y)-u(x)}{|y-x|^{n+s}}\frac{y-x}{|y-x|}\,dy+u(x)c_{n,s}\int_{\Omega^c}\frac{x-y}{|x-y|^{n+s+1}}\,dy\\[1ex]
    &=D^s_\O u(x)+u(x)D^s\mathbbm{1}_\O(x).
\end{align*} 
This identity also extends in a weak sense, that is, for $u \in L^1(\O)$ it holds that $D^s_\O u = D^s(\mathbbm{1}_\O u)-uD^s \mathbbm{1}_\O$ on $\O$, when either $D^s_\O u$ or $D^s(\mathbbm{1}_\O u)|_{\O}$ lies in $L^1_{\rm loc}(\O;\R^n)$.

\subsection{Functions with zero regional fractional gradient}\label{sec:nullspace}
This section is devoted to proving that the regional fractional gradient vanishes exactly for constant functions. The case of regular functions can be proven as follows.
\begin{theorem}
Let $s \in (0,1)$, $\Omega$ be a bounded $C^1$-domain and $u \in C^2(\overline{\Omega})$. Then, it holds that
\[
\text{$D^s_\O u =0$ on $\O$} \quad \text{if and only if} \quad \text{$u$ is constant on $\O$}.
\]
\end{theorem}
\begin{proof}
One implication is clear, so we assume $D^s_\O u=0$ on $\O$. In particular, we find that $\Div D^s_\O u =0$ on $\O$ as well. To compute this divergence, we observe
\[
D^s_\O u(x) = \int_{\O} (u(y)-u(x))\frac{y-x}{\abs{y-x}^{n+s+1}}\,dy=\lim_{r \downarrow 0} \int_{\O \setminus B_r(x)} (u(y)-u(x))\frac{y-x}{\abs{y-x}^{n+s+1}}\,dy,
\]
where we leave out the constant $c_{n,s}$ throughout this proof. We can calculate
\begin{align*}
\Div_x  \int_{\O \setminus B_r(x)}& (u(y)-u(x))\frac{y-x}{\abs{y-x}^{n+s+1}}\,dy \\
&= \int_{\O \setminus B_r(x)} (s+1)\frac{u(y)-u(x)}{\abs{y-x}^{n+s+1}} - \nabla u(x) \cdot \frac{y-x}{\abs{y-x}^{n+s+1}} \,dy \\
&\qquad - \int_{\partial B_r(x)} (u(y)-u(x)) \frac{y-x}{\abs{y-x}^{n+s+1}} \cdot \frac{y-x}{\abs{y-x}}\,d\Hcal^{n-1}(y).
\end{align*}
For the first term, we find
\begin{align*}
\int_{\O \setminus B_r(x)} &(s+1)\frac{u(y)-u(x)}{\abs{y-x}^{n+s+1}} - \nabla u(x) \cdot \frac{y-x}{\abs{y-x}^{n+s+1}} \,dy \\
&=\int_{\O \setminus B_r(x)}(s+1) \frac{u(y)-u(x)-\nabla u(x) \cdot (y-x)}{\abs{y-x}^{n+s+1}} + s\nabla u(x) \cdot \frac{y-x}{\abs{y-x}^{n+s+1}} \,dy.
\end{align*}
As $r \to 0$, using symmetry, this quantity converges uniformly on compact subsets of $\O$ to
\begin{align*}
&(s+1) \int_{\O}  \frac{u(y)-u(x)-\nabla u(x) \cdot (y-x)}{\abs{y-x}^{n+s+1}}\,dy + s\nabla u(x) \cdot \int_{\O \setminus B_{d(x)}(x)}\frac{y-x}{\abs{y-x}^{n+s+1}}\,dy,
\end{align*}
with $d(x):=\dist(x,\O^c)$. For the other term we write
\begin{align*}
\int_{\partial B_r(x)}& (u(y)-u(x)) \frac{y-x}{\abs{y-x}^{n+s+1}} \cdot \frac{y-x}{\abs{y-x}}\,d\Hcal^{n-1}(y)\\
&=\int_{\partial B_r(x)} (u(y)-u(x)-\nabla u(x)\cdot (y-x)) \frac{1}{\abs{y-x}^{n+s}}\,d\Hcal^{n-1}(y),
\end{align*}
which converges uniformly to zero as $r \to 0$. All in all, this proves that
\begin{align}
\begin{split}\label{eq:divdso}
\Div D^s_\O u (x) &=(s+1) \int_{\O}  \frac{u(y)-u(x)-\nabla u(x) \cdot (y-x)}{\abs{y-x}^{n+s+1}}\,dy\\
&\qquad+ s\nabla u(x) \cdot \int_{\O \setminus B_{d(x)}(x)}\frac{y-x}{\abs{y-x}^{n+s+1}}\,dy.
\end{split}
\end{align}
Suppose now that $u$ attains its global minimum at an interior point $x_0 \in \O$. Then, we find that $\nabla u(x_0) = 0$ and thus
\[
\Div D^s_\O u (x_0) =(s+1)\int_{\O} \frac{u(y)-u(x_0)}{\abs{y-x_0}^{n+s+1}}\,dy=0.
\]
Combining this with $u(y) \geq u(x_0)$ for all $y \in \O$, we infer that $u(y)=u(x_0)$ for all $y \in \O$. Hence, $u$ is constant on $\Omega$.

It remains to check the case where the minimum of $u$ is attained at a point $x_0 \in \partial \O$. Since $\partial \O$ is a $C^1$-manifold without boundary, we deduce that $\nabla u(x_0) \cdot v = 0$ for any vector $v \in \R^n$ tangent to $\partial \O$ at $x_0$. Additionally, let $(x_k)_k \subset \O$ be a sequence converging to $x_0$, then we find
\begin{align}
\begin{split}\label{eq:boundary}
0 = \lim_{k \to \infty} \Div D^s_\O u(x_k) &= \lim_{k \to \infty} (s+1)\int_{\O} \frac{u(y)-u(x_k)-\nabla u(x_k) \cdot (y-x_k)}{\abs{y-x_k}^{n+s+1}}\,dy \\
&\qquad+ s\nabla u(x_k) \cdot \int_{\O \setminus B_{d(x_k)}(x_k)}\frac{y-x_k}{\abs{y-x_k}^{n+s+1}}\,dy.
\end{split}
\end{align}
The first term simply converges to
\begin{align*}
&\lim_{k \to \infty} (s+1)\int_{\O} \frac{u(y)-u(x_k)-\nabla u(x_k) \cdot (y-x_k)}{\abs{y-x_k}^{n+s+1}}\,dy \\
&\qquad= (s+1)\int_{\O} \frac{u(y)-u(x_0)-\nabla u(x_0) \cdot (y-x_0)}{\abs{y-x_0}^{n+s+1}}\,dy.
\end{align*}
For the second term, we use the Gauss-Green theorem to find
\begin{align*}
V(x):=&\int_{\O \setminus B_{d(x)}(x)}\frac{y-x}{\abs{y-x}^{n+s+1}}\,dy = \frac{1}{n+s-1}\int_{\partial \O} \frac{\nu_\O(y)}{\abs{y-x}^{n+s-1}}\,d\Hcal^{n-1}(y),
\end{align*}
where $\nu_\Omega$ denotes the inward unit normal to $\partial \O$, and where the surface integral over $\partial B_{d(x)}(x)$ vanishes due to symmetry. We can compute that
\begin{align*}
\abs{V(x_k)\cdot \nu_\O(x_0)} &\geq \frac{1}{n+s-1}\int_{\partial \O}\frac{1}{\abs{y-x_k}^{n+s-1}}\,d\Hcal^{n-1}(y)\\
&\qquad -\frac{1}{n+s-1}\int_{\partial \O}\frac{\abs{\nu_\O(y)\cdot \nu_\O(x_0)-1}}{\abs{y-x_k}^{n+s-1}}\,d\Hcal^{n-1}(y).
\end{align*}
Using continuity of the normal $\nu_\O$ at $x_0$, we infer that $\abs{V(x_k)\cdot \nu_\O(x_0)} \to \infty$ as $k \to \infty$. Hence, in order for \eqref{eq:boundary} to hold we must have that $\lim_{k \to \infty} \nabla u(x_k) \cdot \nu_\O(x_0)=0$, which implies $\nabla u(x_0) \cdot \nu_\O(x_0)=0$. Since $\nabla u(x_0)$ was also zero tangent to $\partial \O$, we conclude that $\nabla u(x_0)=0$. Another computation now yields
\begin{align*}
\abs{V(x_k)} &\leq \frac{1}{n+s-1}\int_{\partial \O}\frac{1}{\abs{y-x_k}^{n+s-1}}\,d\Hcal^{n-1}(y) \\
&\leq \frac{1}{d(x_k)^{1-\epsilon}} \int_{\partial \O}\frac{1}{\abs{y-x_k}^{n+s+\epsilon-2}}\,d\Hcal^{n-1}(y)\\
&\leq \frac{C}{d(x_k)^{1-\epsilon}},
\end{align*}
for $\epsilon \in (0,1-s)$. Hence, the second term in the right-hand side of \eqref{eq:boundary} converges to zero as $k \to \infty$. We thus conclude that 
\[
\int_{\O} \frac{u(y)-u(x_0)}{\abs{y-x_0}^{n+s+1}}\,dy=0,
\]
and now we can argue in the same way as the case of an interior minimum.
\end{proof}

There is the following alternative argument, which does not require $u$ to be smooth, but only works for convex domains. It is inspired by the identity in \cite{CoM15} for the fractional Laplacian of a composition with a convex function, and uses a suitable power of the distance function to the boundary as a test function.
\begin{proposition}\label{prop:nullspace}
    Let $\Omega \subset \R^n$ be a bounded open and convex set. If $u \in L^1(\Omega)$ satisfies $D^s_\O u=0$ on $\O$, then $u$ is constant a.e.~on $\O$.
\end{proposition}
\begin{proof}
First, we note that by \eqref{eq:connDsobarGso} 
\begin{equation}\label{eq:interiorregularity}
\nabla I_{1-s}(\mathbbm{1}_\O u) = D^s(\mathbbm{1}_\O u) = uD^s\mathbbm{1}_\O \quad \text{in $\O$.}
\end{equation}
Since $u \in L^1(\O)$, we deduce by smoothness of $D^s\mathbbm{1}_\O$ that $uD^s\mathbbm{1}_\O \in L^1_{\rm loc}(\O;\R^n) \subset H^{-r}_{\rm loc}(\Omega;\R^n)$ for $r>n/2$. We infer from \eqref{eq:interiorregularity} that $I_{1-s}(\mathbbm{1}_\O u) \in H^{-r+1}_{\rm loc}(\O)$ so that, by the mapping properties of the Riesz potential, we find $u \in H^{-r+s}_{\rm loc}(\O)$, see~e.g.~\cite[Theorem 1.1.25]{Lerner2010}. This, in turn, yields $uD^s\mathbbm{1}_\O \in H^{-r+s}_{\rm loc}(\Omega;\R^n)$, and now we can use a bootstrapping argument to deduce that $u$ must be smooth in the interior, that is, $u \in C^{\infty}(\O)$. Hence, we find that for all $x \in \Omega$ (again leaving out $c_{n,s}$ for ease of notation)
\[
D^s_\O u(x)=\int_{\O}\frac{u(y)-u(x)}{|y-x|^{n+s}}\frac{y-x}{|y-x|}\,dy =0
\]
and
\[
\Div D^s_\O u(x) = (1+s) \lim_{r \downarrow 0} \int_{\O \setminus B_r(x)} \frac{u(y)-u(x)}{|y-x|^{n+s+1}}\,dy - \nabla u(x) \cdot \int_{\Omega^c}\frac{y-x}{|y-x|^{n+s+1}}\,dy =0,
\]
where we have slightly rewritten the expression in \eqref{eq:divdso} using that for each $r>0$
\[
\int_{\R^n \setminus B_{r}(x)} \frac{y-x}{|y-x|^{n+s+1}} \,dy =0.
\]
Now consider the function $F(t):=\sqrt{1+t^2}$ which is strictly convex and has linear growth, which implies that
\[
B(x,y):=F(u(y))-F(u(x))-F'(u(x))(u(y)-u(x)) \geq 0,
\]
and it vanishes if and only if $u(x)=u(y)$. Moreover, with $v:=F \circ u \in L^1(\O) \cap C^\infty(\O)$ we compute
\begin{align*}
D^s_\O v(x) &= \int_{\O} \frac{F(u(y))-F(u(x))}{|y-x|^{n+s}}\frac{y-x}{|y-x|}\,dy\\
&=\int_{\O} \frac{F(u(y))-F(u(x))-F'(u(x))(u(y)-u(x))}{|y-x|^{n+s}}\frac{y-x}{|y-x|}\,dy\\
&=\int_{\O} \frac{B(x,y)(y-x)}{|y-x|^{n+s+1}}\,dy,
\end{align*}
where in the second line we subtracted $F'(u(x))D^s_\O u(x)=0$. For the divergence, we find
\begin{align*}
\Div D^s_\O v(x) &= (1+s)\lim_{r \downarrow 0}\int_{\O \setminus B_r(x)} \frac{F(u(y))-F(u(x))}{|y-x|^{n+s+1}}\,dy\\
&\qquad -F'(u(x))\nabla u(x) \cdot \int_{\Omega^c}\frac{y-x}{|y-x|^{n+s+1}}\,dy \\[1ex]
&= (1+s)\lim_{r \downarrow 0}\int_{\O \setminus B_r(x)} \frac{F(u(y))-F(u(x))-F'(u(x))(u(y)-u(x))}{|y-x|^{n+s+1}}\,dy\\
&= (1+s)\int_{\O} \frac{B(x,y)}{|y-x|^{n+s+1}}\,dy,
\end{align*}
where in the second equality we subtracted $F'(u(x))\Div D^s_\O u(x)=0$ and the last equality uses the non-negativity of $B$ to take the limit $r \downarrow 0$ via the monotone convergence theorem. Next, we consider $d(x):=\dist(x,\O^c)$, which is a concave function on $\O$, and thus, satisfies
\begin{equation}\label{eq:distconc}
\nabla d(x) \cdot (y-x) \geq d(y)-d(x) \quad \text{for $x,y \in \O$.}
\end{equation}
For $\epsilon >0$, we then define the function
\[
\phi_\epsilon(x):=\max\{d(x)-\epsilon,0\}^{1+s},
\]
which is a Lipschitz function with compact support in $\O$. Hence, using integration by parts, we find
\begin{align*}
    (1+s)\int_{\Omega}\int_{\Omega} \frac{B(x,y)}{|y-x|^{n+s+1}}\phi_\epsilon(x)\,dy\,dx & = \int_{\Omega} \Div D^s_\O v(x) \phi_\epsilon(x)\,dx \\
    &=-\int_{\Omega} D^s_\O v(x) \cdot \nabla \phi_\epsilon (x)\,dx \\
    &=-\int_{\O}\int_{\O}\frac{B(x,y)(y-x)}{|y-x|^{n+s+1}}\cdot \nabla \phi_\epsilon(x)\,dy\,dx.
\end{align*}
Subtracting the right-hand side, we find
\begin{align*}
    0&=\int_{\Omega}\int_{\Omega} \frac{B(x,y)}{|y-x|^{n+s+1}}\left((1+s)\phi_\epsilon(x)+\nabla \phi_\epsilon(x)\cdot(y-x)\right)\,dy\,dx \\
    &\geq \int_{\Omega_{-\epsilon}}\int_{\Omega} \frac{B(x,y)}{|y-x|^{n+s+1}}(1+s)(d(x)-\epsilon)^s(d(y)-\epsilon)\,dy\,dx,
\end{align*}
where $\Omega_{-\epsilon}:=\{x \in \O \,:\, d(x)>\epsilon\}$ and we used the chain rule to compute $\nabla \phi_\epsilon$ and applied \eqref{eq:distconc}. We claim that the negative part of this last integral vanishes as $\epsilon \to 0$, that is,
\begin{equation}\label{eq:layervanishing}
\lim_{\epsilon \to 0} \int_{\Omega_{-\epsilon}}\int_{\Omega \setminus \Omega_{-\epsilon}} \frac{B(x,y)}{|y-x|^{n+s+1}}(1+s)(d(x)-\epsilon)^s(\epsilon-d(y))\,dy\,dx=0.
\end{equation}
If this is true, then we may apply Fatou's lemma to the non-negative part and obtain
\[
0 \geq \int_{\Omega}\int_{\Omega} \frac{B(x,y)}{|y-x|^{n+s+1}}(1+s)d(x)^sd(y)\,dy\,dx,
\]
which yields $B(x,y)=0$ for a.e.~$x,y \in \Omega$. In particular, $u(x)=u(y)$ for a.e.~$x,y \in \Omega$ showing that $u$ is constant.

It remains to prove \eqref{eq:layervanishing}. Since $F$ is $1$-Lipschitz, we find that
\[
0 \leq B(x,y) \leq 2|u(y)-u(x)| \leq 2 \left(|u(x)|+|u(y)|\right).
\]
Additionally, for every $x \in \Omega_{-\epsilon}$ we can estimate
\begin{align*}
    \int_{\Omega \setminus \Omega_{-\epsilon}} \frac{(d(x)-\epsilon)^s(\epsilon-d(y))}{|y-x|^{n+s+1}}\,dy & \leq \int_{\Omega \setminus \Omega_{-\epsilon}} \frac{(d(x)-\epsilon)^s\min\{|y-x|,\epsilon\}}{|y-x|^{n+s+1}}\,dy \\
    & \leq C (d(x)-\epsilon)^s \min\{(d(x)-\epsilon)^{-s},\epsilon (d(x)-\epsilon)^{-(1+s)}\} \\
    &= C \min\{1,\epsilon (d(x)-\epsilon)^{-1}\},
\end{align*}
where we used 
\[
|y-x| \geq d(x)-d(y) \geq \max\{d(x)-\epsilon,\epsilon-d(y)\} \quad \text{for all $x \in \Omega_{-\epsilon}$ and $y \in \Omega \setminus \Omega_{-\epsilon}$.}
\]
This shows via Lebesgue's dominated convergence theorem that
\[
\lim_{\epsilon\to 0} \int_{\Omega_{-\epsilon}} |u(x)| \int_{\Omega \setminus \Omega_{-\epsilon}} \frac{(d(x)-\epsilon)^s(\epsilon-d(y))}{|y-x|^{n+s+1}}\,dy\,dx=0.
\]
For the term with $|u(y)|$, we reverse the integrals via Fubini's theorem and estimate in a similar way
\begin{align*}
    \int_{\Omega_{-\epsilon}}\frac{(d(x)-\epsilon)^s(\epsilon-d(y))}{|y-x|^{n+s+1}}\,dx &\leq \int_{\Omega_{-\epsilon}}\frac{(\epsilon-d(y))}{|y-x|^{n+1}}\,dx \\
    &\leq C (\epsilon-d(y))^{-1}(\epsilon-d(y))\leq C.
\end{align*}
Hence, 
\[
\lim_{\epsilon\to 0} \int_{\Omega \setminus \Omega_{-\epsilon}} |u(y)| \int_{\Omega_{-\epsilon}} \frac{(d(x)-\epsilon)^s(\epsilon-d(y))}{|y-x|^{n+s+1}}\,dx\,dy \leq C\lim_{\epsilon\to 0} \int_{\Omega \setminus \Omega_{-\epsilon}} |u(y)|\,dy=0.
\]
This finishes the proof of \eqref{eq:layervanishing} and the result.
\end{proof}

A different observation complementing the previous results is the following.
\begin{proposition}
    Let $\Omega$ be a bounded open set and $u\in L^1(\O)$. If $u$ is constant a.e. in a non-empty open subset $\Omega_0\subset\Omega$, and $D^s_\Omega u=0$ in $\Omega_0$, then $u$ is constant a.e. in $\Omega$.
\end{proposition}
\begin{proof}
    Let $c \in \R$ be the constant value of $u$ in $\O_0$ and $v$ be the extension by zero of $u-c$. Hence, $v\in L^1(\R^n)$, which implies that $v\in H^{-r}(\R^n)$ for every $r>n/2$. Now, we have that $$D^sv=D^s_\Omega(u-c)+(u-c)D^s\mathbbm{1}_\Omega,$$ and if we restrict to $x\in \Omega_0$, we get that $D^sv=0$ in $\Omega_0$. Now, using the identity $\operatorname{div}D^s=-(-\Delta)^{(1+s)/2}$, we get that both $v$ and $(-\Delta)^{(1+s)/2}v$ are zero in $\Omega_0$. This implies by \cite[Theorem 1.2]{ghosh2020}, that $v=0$ in the whole space, so, in particular, $u=c$ in the whole domain $\Omega$.
\end{proof}

\subsection{Bessel potential spaces on domains and the regional spaces}\label{sec:connection}
This section is devoted to the connections between the Bessel potential spaces on domains and the regional fractional spaces. The following result shows that the Bessel potential spaces on domains are embedded into the regional fractional spaces.
\begin{theorem}\label{th:emb1}
    Let $s \in (0,1)$, $p \in (1,\infty)$ and $\Omega$ be an open and bounded set. Then, $$H^{s,p}(\Omega)\xhookrightarrow{}\mathcal{H}^{s,p}(\Omega).$$
\end{theorem}
\begin{proof}
    We first recall the result from \cite[Corollary 1.4]{li2019}, which states that given a differential operator $A^s$ such that its symbol is homogeneous of degree $s>0$ and $\widehat{A}^s(\xi)\in C^\infty(\mathbb{S}^{n-1})$, as, for example, the symbol $2\pi i \xi_j|2\pi \xi|^{s-1}$ of the $j$th component of the Riesz fractional gradient, then it holds that
    \[
    \norm{A^s(fg)-fA^s(g)-gA^s(f)}_p\lesssim \norm{\ds f}_p\norm{g}_{BMO} \quad \text{for all $f,g \in \Scal(\R^n)$.}
    \]
    By the triangle inequality, 
    \[
    \norm{A^s(fg)-fA^s(g)}_p\lesssim \norm{\ds f}_p\norm{g}_{BMO}+\norm{g}_\infty \norm{A^s(f)}_p.
    \]
    Since for functions in $L^\infty(\R^n)$ we have that $\norm{g}_{BMO}\leq 2\norm{g}_\infty$, and by choosing $A^s=D^s$ and using $D^sf=\mathcal{R}\ds f$ together with the boundedness of the Riesz transform, we get that 
    \[
    \norm{D^s (fg)-fD^sg}_p\lesssim \norm{\ds f}_p\norm{g}_\infty.
    \]
    Now, let $(g_j)_j \subset C_c^{\infty}(\R^n)$ be a sequence that converges to $\mathbbm{1}_\Omega$ pointwise a.e.~and satisfies $\|g_j\|_\infty \leq 1$ for all $j \in \N$ and $g_j \equiv 1$ on $\Omega_j:=\{x \in \Omega\,:\, \dist(x,\Omega^c) \geq 1/j\}$, which, for example, can be constructed via mollification. Then, we find that
    \[
    \norm{D^s (fg_j)-fD^sg_j}_p\lesssim \norm{\ds f}_p\norm{g_j}_\infty \lesssim \norm{\ds f}_p.
    \]
    In particular, we find that $V_j:=D^s (fg_j)-fD^sg_j$ converges weakly, up to a subsequence, to a $V \in L^p(\R^n;\R^n)$ with
    \[
    \|V\|_p \leq \liminf_{j \to \infty} \|V_j\|_p \lesssim \norm{\ds f}_p.
    \]
    Next, we show using a distributional argument that $V|_{\Omega}=D^s_\Omega (f|_{\Omega})$. To this aim, we first show that for any compact set $K \Subset \Omega$ we have that
    \begin{equation}\label{eq:indicatorconvergence}
    \|D^s g_j \|_{L^{\infty}(K)}\leq C \quad \text{and} \quad D^s g_j \to D^s \mathbbm{1}_\Omega \ \text{pointwise on $K$.}
    \end{equation}
    Indeed, for $j \geq k_0> \dist(K,\Omega^c)^{-1}$ and $x \in K$ we infer that
    \[
    |D^s g_j (x)| = c_{n,s}\abslr{\int_{\Omega_{k_0}^c} \frac{(1-g_j(y))(y-x)}{|y-x|^{n+s+1}}\,dy} \leq 2c_{n,s}\int_{\Omega_{k_0}^c} \frac{1}{|y-x|^{n+s}}\,dy < \infty,
    \]
which shows the uniform bound and allows us to conclude the pointwise convergence using Lebesgue's dominated convergence theorem.
Now we take $\phi \in C_c^{\infty}(\Omega;\R^n)$, then we can compute that
    \begin{align*}
        \int_{\Omega} V \cdot \phi \,dx &= \lim_{j \to \infty} \int_{\Omega} \left(D^s (fg_j)-fD^sg_j \right) \cdot \phi \,dx\\
        &=\lim_{j \to \infty}\int_{\R^n} -fg_j  \Div^s \phi -fD^sg_j \cdot \phi \,dx
        \\
        &=\int_{\Omega} -f  \Div^s \phi -fD^s \mathbbm{1}_\Omega \cdot \phi \,dx\\
        &=- \int_{\Omega} f \Div^s_\Omega \phi \,dx,
    \end{align*}
where in the third line we use Lebesgue's dominated convergence theorem together with the pointwise convergence $g_j \to \mathbbm{1}_\Omega$ and \eqref{eq:indicatorconvergence}. We conclude that $V|_{\Omega} = D^s_\Omega (f|_{\Omega})$, and hence,
\begin{equation}\label{eq:smoothbound}
\|D^s_\Omega (f|_{\Omega})\|_{L^p(\Omega)} \leq \|V\|_p \lesssim \|\ds f\|_p \lesssim \|f\|_{H^{s,p}(\R^n)} \quad \text{for all $f \in \Scal(\R^n)$.}
\end{equation}

To finish the proof, we consider now $u \in H^{s,p}(\Omega)$ and an extension $U \in H^{s,p}(\R^n)$ with $\|U\|_{H^{s,p}(\R^n)}\lesssim \|u\|_{H^{s,p}(\Omega)}$. Moreover, we consider a sequence $(U_j)_j \subset \Scal(\R^n)$ with $U_j \to U$ in $H^{s,p}(\R^n)$. Then, we find using \eqref{eq:smoothbound} that
\[
\|D^s_\Omega(U_j|_{\Omega})-D^s_\Omega(U_k|_{\Omega})\|_{L^p(\Omega)} \lesssim \|U_j-U_k\|_{H^{s,p}(\R^n)} \quad \text{for $j,k \in \N$,}
\]
which shows that $(U_j)_j$ is a Cauchy sequence in $\Hcal^{s,p}(\Omega)$. By completeness of $\Hcal^{s,p}(\O)$ (Proposition~\ref{prop:banach}), we find that $U_j \to u \in \Hcal^{s,p}(\Omega)$ and
\[
\|D^s_\Omega u\|_{L^p(\Omega)} = \lim_{j \to \infty} \|D^s_\Omega(U_j|_{\Omega})\|_{L^p(\Omega)} \leq \lim_{j \to \infty} \|U_j\|_{H^{s,p}(\R^n)} \lesssim \|u\|_{H^{s,p}(\Omega)},
\]
which finishes the proof.
\end{proof}

A specific consequence of the previous result is the following type of fractional Hardy inequality for functions with non-zero trace. To the best of the authors' knowledge, this inequality has never appeared in the literature, yet it shares some similarities with the trace conjunction inequalities recently established in \cite{Sch25}.
\begin{corollary}\label{cor:hardy}
    Let $s \in (0,1)$, $p \in (1,\infty)$ with $sp>1$ and suppose that $\Omega \subset \R^n$ is a bounded Lipschitz domain. Then, there exists a constant $C>0$ such that
    \begin{equation}\label{eq:hardynonzero}
    \int_\O \abslr{\int_{\partial \O} \nu_\O(y)\frac{\mathrm{Tr}\,u(y)-u(x)}{|y-x|^{n+s-1}}\,d\Hcal^{n-1}(y)}^p\,dx \leq C \|u\|_{H^{s,p}(\O)}^p \quad \text{for all $u \in H^{s,p}(\O)$,}
    \end{equation}
    where $\nu_\O$ is the inner unit normal to $\partial\O$ and $\mathrm{Tr}\,u \in W^{s-1/p,p}(\partial\O)$ denotes the trace of $u$.
\end{corollary}
\begin{proof}
Let $u \in C^{\infty}(\overline{\Omega})$. Then, we can compute that
\[
D^s(\mathbbm{1}_{\Omega}u) = I_{1-s}(D(\mathbbm{1}_{\Omega}u)) = I_{1-s}(\mathbbm{1}_{\Omega} \nabla u + \mathrm{Tr}(u)\nu_{\Omega}\,\Hcal^{n-1}\mres_{\partial\Omega}),
\]
with the first identity using Proposition~\ref{DsProps}\,(ii) extended to the $BV$-function $\mathbbm{1}_\O u$, with $D$ denoting its distributional derivative. This shows together with \eqref{eq:connDsobarGso} that
\begin{equation}\label{eq:connDsoGso2}
D^s_\Omega u = I_{1-s}(\mathbbm{1}_{\Omega} \nabla u) +  I_{1-s}(\mathrm{Tr}(u)\nu_{\Omega}\,\Hcal^{n-1}\mres_{\partial\Omega})-uD^s \mathbbm{1}_\Omega.
\end{equation}
For the first term on the right-hand side, we can estimate
\begin{align}\label{eq:estGso}
    \|I_{1-s}(\mathbbm{1}_\Omega \nabla u)\|_{L^p(\O)} \lesssim \|\mathbbm{1}_\O \nabla u\|_{H^{s-1,p}(\R^n)}\lesssim \|\nabla u\|_{H^{s-1,p}(\O)} \lesssim \|u\|_{H^{s,p}(\O)},
\end{align}
where the first inequality uses the mapping properties of the Riesz potential (together with the compact support of $\mathbbm{1}_\O \nabla u$), and the second uses that extension by zero is bounded on $H^{s-1,p}(\O)$ given that $\mathbbm{1}_{\Omega}$ is a multiplier on its dual space $H^{1-s,p'}(\R^n)$ in light of $(1-s)p'<1$. Combining this inequality with
\[
\|D^s_\O u\|_{L^p(\O)} \lesssim \|u\|_{H^{s,p}(\O)},
\]
from Theorem~\ref{th:emb1}, we find that
\[
\|I_{1-s}(\mathrm{Tr}(u)\nu_{\Omega}\,\Hcal^{n-1}\mres_{\partial\Omega})-uD^s \mathbbm{1}_\Omega\|_{L^p(\O)} \lesssim \|u\|_{H^{s,p}(\O)} \quad \text{for all $u \in C^{\infty}(\overline{\O})$.}
\]
Since 
\begin{equation}\label{eq:connDsoGso3}
D^s \mathbbm{1}_\O = I_{1-s}(D \mathbbm{1}_\O)=I_{1-s}(\nu_{\Omega}\,d\Hcal^{n-1}\mres_{\partial\Omega})
\end{equation}
this inequality is exactly \eqref{eq:hardynonzero} up to a constant. By density of $C^{\infty}(\overline{\O})$ in $H^{s,p}(\O)$, this finishes the proof.
\end{proof}
We now give a partial converse to Theorem \ref{th:emb1}, showing that functions in the regional fractional space have local Bessel regularity.

\begin{theorem}\label{th:regional-local-bessel}
Let $s\in(0,1)$, $p\in(1,\infty)$ and $\Omega\subset\mathbb R^n$ be an open and bounded set. Then,
\[
\mathcal H^{s,p}(\Omega)
\hookrightarrow
L^p(\Omega)\cap H^{s,p}_{\rm loc}(\Omega).
\]
More precisely, for every $\eta\in C_c^\infty(\Omega)$ there exists
a constant $C>0$ such that
\[
\|\eta u\|_{H^{s,p}(\mathbb R^n)}
\leq C\left(
\|u\|_{L^p(\Omega)}
+\|D^s_\Omega u\|_{L^p(\Omega;\mathbb R^n)}
\right)
\quad\text{for all }u\in\mathcal H^{s,p}(\Omega),
\]
where $\eta u$ is extended by zero outside $\Omega$.
\end{theorem}

\begin{proof}
Let $u\in\mathcal H^{s,p}(\Omega)$ and $\bar{u}$ be its extension to $\R^n$ by zero, then we use the identity \eqref{eq:connDsobarGso} to find
\[
D^s \bar{u}
=
D^s_\Omega u+uD^s\mathbbm{1}_\Omega
\quad\text{in }\Omega.
\]
For $x\in\Omega$, we can estimate
\begin{equation}\label{eq:distest}
\begin{aligned}
|D^s\mathbbm{1}_\Omega(x)|
&=
c_{n,s}\left|
\int_{\Omega^c}\frac{x-y}{|x-y|^{n+s+1}}\,dy
\right|\\
&\leq
c_{n,s}\int_{|z|\geq d(x)}|z|^{-n-s}\,dz
=
\frac{c_{n,s}|\mathbb S^{n-1}|}{s}\,d(x)^{-s}.
\end{aligned}
\end{equation}
Consequently, $D^s\mathbbm{1}_\Omega$ is bounded on every compact
subset of $\Omega$, and hence
\[
D^s \bar{u}|_{\O} \in L^p_{\mathrm{loc}}(\Omega;\mathbb R^n),
\]
which implies that for any fixed $\eta\in C_c^\infty(\Omega)$, extended by zero to
$\mathbb R^n$, there exists a positive constant $C$ such that 
\[
\|\eta D^s \bar{u}\|_{L^p(\mathbb R^n;\mathbb R^n)}
\leq C\left(
\|D^s_\Omega u\|_{L^p(\Omega;\mathbb R^n)}
+\|u\|_{L^p(\Omega)}
\right).
\]
Now, we observe that
\[
D^s(\eta\bar u)(x)
=
\eta(x)D^s\bar u(x)
+
c_{n,s}\int_{\mathbb R^n}
\bigl(\eta(y)-\eta(x)\bigr)\bar u(y)
\frac{y-x}{|y-x|^{n+s+1}}\,dy.
\]
This identity follows from the integration by parts formula;
the integrable bound below justifies the change in the order
of integration.

Indeed, since $\eta$ is bounded and Lipschitz,
\[
|\eta(y)-\eta(x)|
\leq C_\eta\min\{|y-x|,1\},
\]
therefore,
\[
\left|
\int_{\mathbb R^n}
\bigl(\eta(y)-\eta(x)\bigr)\bar u(y)
\frac{y-x}{|y-x|^{n+s+1}}\,dy
\right|
\leq C_\eta(k*|\bar u|)(x),
\]
where
\[
k(z):=\frac{\min\{|z|,1\}}{|z|^{n+s}}.
\]
Since $s \in (0,1)$, $k$ is integrable and hence 
Young's inequality yields
\[
\|D^s(\eta\bar u)\|_{L^p(\mathbb R^n;\mathbb R^n)}
\leq C\left(
\|D^s_\Omega u\|_{L^p(\Omega;\mathbb R^n)}
+\|u\|_{L^p(\Omega)}
\right).
\]
Since also $\eta\bar u\in L^p(\mathbb R^n)$, we know that
$\eta\bar u\in H^{s,p}(\mathbb R^n)$ so
\[
\begin{aligned}
\|\eta\bar u\|_{H^{s,p}(\mathbb R^n)}
&\leq C\left(
\|\eta\bar u\|_{L^p(\mathbb R^n)}
+\|D^s(\eta\bar u)\|_{L^p(\mathbb R^n;\mathbb R^n)}
\right)\\
&\leq C\left(
\|u\|_{L^p(\Omega)}
+\|D^s_\Omega u\|_{L^p(\Omega;\mathbb R^n)}
\right).
\end{aligned}
\]
As $\eta\in C_c^\infty(\Omega)$ was arbitrary, this proves
$u\in H^{s,p}_{\mathrm{loc}}(\Omega)$.
\end{proof}

In general $H^{s,p}(\O)\subset L^p(\Omega)\cap H^{s,p}_{\rm loc}(\Omega)$ with strict inclusion, as the following example shows.
\begin{example}
For $s\in (0,1)$ and $p \in (1,\infty)$, we define the function $u(x)=x^{-\alpha}$ for $x\in (0,1)$ and for some
    \[
    \max\left\{0,\frac1p-s\right\}<\alpha<\frac1p.
    \]
    Then, we can prove that
    \[
    u\in L^p(0,1)\cap H^{s,p}_{\mathrm{loc}}(0,1),
    \quad \text{and} \quad
    u\notin H^{s,p}(0,1).
    \]
    Indeed, since $\alpha p<1$,
    \[
    \int_0^1|u(x)|^p\,dx
    =\int_0^1x^{-\alpha p}\,dx
    =\frac1{1-\alpha p}<\infty,
    \]
    so $u\in L^p(0,1)$. Also, $u$ is smooth inside $(0,1)$, so $u\in H^{s,p}_{\rm loc}(0,1)$.
    Since
    \[
    0<\frac1p-\alpha<\min\left\{s,\frac1p\right\},
    \]
    we can choose $t$ such that
    \[
    \frac1p-\alpha<t<\min\left\{s,\frac1p\right\},
    \]
    and $q$ such that
    \[
    \frac1q=\frac1p-t,
    \qquad\text{i.e.,}\qquad
    q=\frac{p}{1-tp},
    \]
    satisfying $0<t<s<1$, $tp<1$ and $p<q<\infty$.

    Since $u\in H^{s,p}(0,1)$ implies $u\in H^{t,p}(0,1)$ and since $tp<1$, the subcritical Sobolev embedding gives us that $u\in L^q(0,1)$. But the choice of $t$ implies that
    \[
    \frac1q=\frac1p-t<\alpha
    \quad\Rightarrow \quad \alpha q>1.
    \]
    As a consequence,
    \[
    \int_0^1|u(x)|^q\,dx
    =\int_0^1x^{-\alpha q}\,dx=\infty,
    \]
    which is a contradiction. We conclude that $u\notin H^{s,p}(0,1)$.
\end{example}

\begin{remark}[Truncated fractional gradients]
We recall the relationship between the regional fractional space $\mathcal{H}^{s,p}(\Omega)$ with the spaces $H^{s,p,\delta}(\Omega)$ introduced in \cite{BeCuMC23} as the space associated to the \textit{truncated Riesz fractional gradient} $D^s_\delta$. This nonlocal operator was introduced as an alternative to the Riesz fractional gradient $D^s$ suitable for bounded domains. It consists of a truncation of $D^s$ via a smooth cut-off function which reduces the interaction between points to a finite horizon $\delta>0$. This allows us to work on bounded domains $\Omega$, albeit with the disadvantage of requiring the function $u$ to be defined in a collar of width $\delta$ around the domain in order for $D^s_\delta u$ to be well-defined in $\Omega$. Since the spaces shrink as $\delta$ approaches zero, we can consider the intersection of all of them, and from there it is straightforward to prove
\[
\mathcal{H}^{s,p}(\Omega)\xhookrightarrow{}\bigcap_{\delta>0}H^{s,p,\delta}(\Omega_{-\delta}) \quad \text{and} \quad \bigcap_{\delta>0}H^{s,p,\delta}(\Omega_{-\delta})=L^p(\Omega)\cap H^{s,p}_{\rm loc}(\Omega),
\]
with $\Omega_{-\delta}:=\{x \in \O \,:\, d(x,\partial \O) > \delta\}$. We refer to \cite{GuilleThesis} for more details on this.
\end{remark}

The full identification with the Bessel potential spaces holds up to additional control on the growth of the function near the boundary in the form of a fractional Hardy inequality. This yields the following novel characterization of Bessel potential spaces on domains. The proof is postponed to the end of Section~\ref{sec:extended}, since the techniques rely on the results for the extended fractional gradients that are developed there.

\begin{theorem}[Characterization of Bessel potential spaces]\label{th:char}
    Let $s \in (0,1)$ and $p \in (1,\infty)$ with $p\not = 2/(1+s)$, and $\Omega \subset \R^n$ be a smooth and bounded domain. Then, the following holds:
    \begin{itemize}
        \item[(i)] It holds that $u \in H^{s,p}_0(\O)$ if and only if
        \[
        u \in \Hcal^{s,p}(\O) \quad \text{and} \quad u/d^s \in L^p(\O).
        \]
        Moreover, there is a $C>0$ such that
        \[
        \|u\|_{H^{s,p}(\R^n)} \leq C\left(\|u/d^s\|_{L^p(\O)}+\|D^s_\O u\|_{L^p(\O;\R^n)}\right) \quad\text{for all $u \in H^{s,p}_0(\O)$.}
        \]

        \item[(ii)] If $sp>1$, then $u \in H^{s,p}(\O)$ if and only if
        \[
        u \in \Hcal^{s,p}(\O) \quad \text{and} \quad \int_{\partial \Omega}\nu_\Omega(y) \frac{g(y)-u(\cdot)}{|y-\cdot|^{n+s-1}}\,d\Hcal^{n-1}(y) \in L^p(\O;\R^n),
        \]
        for some $g \in W^{s-1/p,p}(\partial \O)$; if this holds, $g$ is the trace of $u$. Moreover, there is a $C>0$ such that for all $u \in H^{s,p}(\O)$
        \begin{align*}
        &\|u-(u)_\O\|_{H^{s,p}(\O)} \\
        &\qquad \leq C\left(\normlr{\int_{\partial \O} \nu_\O(y)\frac{\mathrm{Tr}\,u(y)-u(\cdot)}{|y-\cdot|^{n+s-1}}\,d\Hcal^{n-1}(y)}_{L^p(\O;\R^n)}+\|D^s_\O u\|_{L^p(\O;\R^n)}\right).
        \end{align*}
    \end{itemize}
\end{theorem}

\section{Characterizations of Bessel potential spaces in domains using extended fractional gradients}\label{sec:extended}
In this section, we characterize Bessel potential spaces on domains by using fractional gradients defined on suitably extended functions. Depending on the values of $s$ and $p$, these distributional spaces characterize the Bessel spaces in domains. Moreover, for almost all parameter values, they provide equivalent norms on the Bessel spaces, leading to new fractional Poincar\'e-Wirtinger inequalities.

The idea of the extended fractional gradients is to suitably combine the full-space fractional operators and zero-extension of functions. Precisely, we consider for $u \in C^{\infty}(\overline{\O})$ the operators
\[
\bar{G}^s_\O u:=D^s(\mathbbm{1}_\O u)|_{\O} \quad \text{and} \quad G^s_\O u:= I_{1-s}(\mathbbm{1}_\O \nabla u)|_{\O}.
\]
The former applies the full-space fractional gradient to the zero-extension, and thereby, sees a possible jump at the boundary. Hence, this gradient is only suitable in the case $sp \leq 1$, when traces do not exist, or for $sp > 1$ and functions with zero trace. The second extended fractional gradient does not penalize such jumps at the boundary, through extension of the gradient instead of the function itself, and thus makes it amenable to functions with non-zero trace. However, for functions in the fractional Sobolev space the gradient only exists distributionally, so its zero-extension is only sensible when $sp>1$, i.e., when traces exist. This is somewhat natural, since the connection between the two gradients is
\begin{equation}\label{eq:identityG}
G^s_\O u = \bar{G}^s_\O u-I_{1-s}(u \nu_\O \,\Hcal^{n-1}\mres_{\partial\O}) \quad \text{for $u \in C^{\infty}(\overline{\O})$,}
\end{equation}
which shows that the second gradient is a correction of the first with a term given by the boundary values of $u$.

We first define the two gradients in a distributional way.
\begin{definition}[Weak extended fractional gradient]
    Let $s \in (0,1)$, $p \in (1,\infty)$, $\Omega \subset \R^n$ be a bounded Lipschitz domain and $u \in L^p(\O)$.
    \begin{itemize}
        \item[(i)] We say that $U=\bar{G}^s_\O u \in L^p(\O;\R^n)$ if
        \[
        \int_{\O} U  \cdot \varphi \,dx = - \int_{\O} u \Div^s \varphi \,dx \quad \text{for all $\varphi \in C_c^{\infty}(\O;\R^n)$.}
        \]

        \item[(ii)] If $sp>1$, we say that $U=G^s_\O u \in L^p(\O;\R^n)$ if there exists a $g \in L^p(\partial\O)$ such that
            \[
    \int_{\O}U \cdot \varphi \,dx = -\int_{\O}u\Div^s \varphi \,dx + \int_{\partial\O}gI_{1-s}\varphi \cdot \nu_\O \,d\Hcal^{n-1} \quad \text{for all $\varphi \in C_c^\infty(\O;\R^n).$}
    \]
    We also write in this case $\Trt u:=g$.
    \end{itemize}
\end{definition}

\begin{definition}[Fractional Sobolev space]
    Let $s \in (0,1)$, $p \in (1,\infty)$ and $\Omega \subset \R^n$ be a bounded Lipschitz domain. We define the space
    \[
    \widetilde{H}^{s,p}_0(\O):=\left\{u \in L^p(\Omega) \,:\, \bar{G}^s_\O u \in L^p(\Omega;\R^n) \right\}
    \]
    and, for $sp>1$,
    \[
    \widetilde{H}^{s,p}(\O):=\{u \in L^p(\O)\,:\, G^s_\O u \in L^p(\O;\R^n)\}.
    \]
\end{definition}

\begin{remark}
    a) We remark that $G^s_\O u$ and $\Trt u$ are actually uniquely defined for $sp>1$. Indeed, by Fubini's theorem, it holds that
    \[
    \int_{\O}(G^s_\O u+I_{1-s}(g \nu_\O \Hcal^{n-1}\mres_{\partial\O})) \cdot \varphi \,dx = -\int_{\O}u\Div^s \varphi \,dx \quad \text{for all $\varphi \in C_c^\infty(\O;\R^n).$}
    \]
    Additionally, $W:=I_{1-s}(g \nu_\O \Hcal^{n-1}\mres_{\partial\O}) \in L^p(\O;\R^n)$ if and only if $g \equiv 0$, given the fact that 
    \[
    \lim_{t \downarrow 0} t^sW(x_0+t\nu_\O(x_0)) \propto g(x_0)\nu_\O(x_0) \quad \text{for every Lebesgue point $x_0 \in \partial\O$ of $g$.}
    \]
    Indeed, this can be proven by using that the measures 
    \[
    \mu_{t,x_0}:=\frac{t^s}{|x_0+t\nu_\O(x_0)-\cdot|^{n+s-1}}\Hcal^{n-1}\mres_{\partial\O}
    \]
    localize to a multiple of the dirac measure at $x_0$ as $t \downarrow 0$. Then, since $sp>1$, we find that $W$ is not $p$-integrable on $\{x_0+t\nu_\O(x_0)\,:\, t >0\}\cap \O$ for $\Hcal^{n-1}$-a.e.~$x_0 \in \partial\O$, which implies that $W$ is not in $L^p(\O;\R^n)$. Therefore, $g=\Trt u$ is uniquely determined by the condition that $G^s_\O u \in L^p(\O;\R^n)$, after which the uniqueness of $G^s_\O u$ also follows. \medskip

    b) It can be readily checked from the definition that if $sp>1$, then
    \[
    \widetilde{H}^{s,p}_0(\O)=\{u \in \widetilde{H}^{s,p}(\O) \,:\, \Trt u =0\} \quad \text{and} \quad \bar{G}^s_\O = G^s_\O \quad \text{on $\widetilde{H}^{s,p}_0(\O)$.} 
    \]
    This is in line with the identity from \eqref{eq:identityG}.
\end{remark}

We have the following characterization using the regularity theory from \cite{Gru15}, where we recall the exponents
\[
q_s=\frac{2}{1+s} \quad \text{and} \quad q_s'=\frac{2}{1-s}.
\]

\begin{theorem}[Characterization of Bessel spaces in domains]\label{th:tildechar}
    Let $s \in (0,1)$, $p \in (1,\infty)$ and $\Omega \subset \R^n$ be a smooth and bounded domain. 
    \begin{itemize}
        \item[(i)] It holds that
    \begin{equation}\label{eq:tildechar-zero}
    \widetilde{H}^{s,p}_0(\Omega) = \left\{u \in L^p(\Omega) \,:\, \mathbbm{1}_\Omega u \in H^{\frac{s-1}{2}(s)}_p(\overline{\Omega}) \right\}.
    \end{equation}
    Moreover, we have that 
    \begin{equation}\label{eq:idnohardy}
    \widetilde{H}^{s,p}_0(\Omega) = H^{s,p}_0(\Omega) \quad \text{when $p < q_s$ or $p> q_s'$.}
    \end{equation}
    Finally, if $u \in \widetilde{H}^{s,p}_0(\Omega)$ satisfies $u/d^s \in L^p(\Omega)$, then also $u \in H^{s,p}_0(\Omega)$ for all $s,p$ with $p \not = q_s$. \smallskip

    \item[(ii)] If $sp>1$, it holds that
    \begin{equation}\label{eq:tildechar-trace}
    \widetilde{H}^{s,p}(\Omega)=H^{s,p}(\Omega) \quad \text{when $p> q_s'$.}
    \end{equation}
    Additionally, if $u \in \widetilde{H}^{s,p}(\O)$ with $\widetilde{\mathrm{Tr}}\,u \in W^{s-1/p,p}(\partial \O)$ and
    \begin{equation}\label{eq:hardytilde}
    \int_{\partial \Omega}\nu_\Omega(y) \frac{\widetilde{\mathrm{Tr}}\,u(y)-u(\cdot)}{|y-\cdot|^{n+s-1}}\,d\Hcal^{n-1}(y) \in L^p(\O;\R^n),
    \end{equation}
    then $u \in H^{s,p}(\O)$ also when $1/s < p \leq q_s'$.
    \end{itemize}
\end{theorem}
\begin{proof}
    \textit{Part $(i)$.} From Theorem~\ref{th:regularity} and \cite[Theorem~4.2]{Gru15} it follows that $I_{1-s}(\mathbbm{1}_{\Omega}u)|_{\Omega} \in W^{1,p}(\Omega)$, that is, $u \in \widetilde{H}^{s,p}_0(\Omega)$, if and only if
    \[
 \mathbbm{1}_\Omega u \in H^{\frac{s-1}{2}(s)}_p(\overline{\Omega}).
    \]
    This shows the identification \eqref{eq:tildechar-zero}.

    To prove \eqref{eq:idnohardy}, we first note that for $u \in H^{s,p}_0(\O)$ we have $\mathbbm{1}_{\O}u \in H^{s,p}(\R^n)$, which shows
    \[
    \bar{G}^s_\O u = D^s(\mathbbm{1}_\O u)|_{\O} \in L^p(\O;\R^n),
    \]
    that is, $u \in \widetilde{H}^{s,p}_0(\O)$. For the reverse inclusion, if $p < q_s$, then $H^{\frac{s-1}{2}(s)}_p(\overline{\Omega})=H^{s,p}_0(\Omega)$, so that \eqref{eq:idnohardy} follows directly from \eqref{eq:tildechar-zero}. For the case $p>q_s'$, we note that $\mathbbm{1}_{\Omega} u \in L^p_0(\Omega) \cap H^{\frac{s-1}{2}(s)}_p(\overline{\Omega})$ implies by \eqref{eq:orderreducing} and Definition~\ref{def:transmission} that
    \[
    \Lambda^{(\frac{s-1}{2})}_{+}(\mathbbm{1}_{\Omega} u) \in H^{\frac{1-s}{2},p}_0(\Omega) \cap H^{\frac{1+s}{2},p}(\Omega)=H^{\frac{1+s}{2},p}_0(\Omega),
    \]
    where the second equality uses that $\frac{1-s}{2}p>1$, so that traces exist. We then conclude from \eqref{eq:orderreducing} that $\mathbbm{1}_\O u \in H^{s,p}_0(\O)$, which proves \eqref{eq:idnohardy}.
    
    Finally, when $p \in (q_s,q_s']$, then any $u \in H^{\frac{s-1}{2}(s)}_p(\overline{\Omega})$ can be decomposed as $u=v+w$ with $\mathbbm{1}_\O w \in H^{s,p}_0(\Omega)$ and $\tilde{v}:=d^{\frac{1-s}{2}}v \in H^{\frac{1+s}{2},p}(\Omega)$ \cite[Theorem~3.4]{Gru19}; in fact, if the trace of $\tilde{v}$ is zero (note that the trace exists since $\frac{1+s}{2}p>1$), we may assume that $v$ is zero as well, which would yield $u = w \in H^{s,p}_0(\Omega)$ as desired. To prove that the trace of $\tilde{v}$ is zero, we note that $\tilde{v}$ satisfies the fractional Hardy inequality; indeed, we have that
    \[
    \frac{\tilde{v}}{d^{\frac{1+s}{2}}} = \frac{v}{d^s}=\frac{u-w}{d^s} \in L^p(\Omega).
    \]
    Hence, we have that $\tilde{v} \in H^{\frac{1+s}{2},p}_0(\Omega)$, which finishes the proof of this case. \medskip

    \textit{Part $(ii)$.} We note first that for $u \in H^{s,p}(\O)$ and $sp>1$, it holds that $I_{1-s}(\mathbbm{1}_\O \nabla u) \in L^p(\O;\R^n)$ by \eqref{eq:estGso}. We want to show that $G^s_\O u =I_{1-s}(\mathbbm{1}_\O \nabla u)$. Since for $u \in C^{\infty}(\overline{\Omega})$ we have
\[
\int_{\R^n} \mathbbm{1}_\Omega \nabla u \cdot \psi \,dx = -\int_{\O}u \Div \psi \,dx +\int_{\partial \O} u \psi \cdot \nu_\O \,d\Hcal^{n-1} \quad \text{for all $\psi \in C^\infty(\R^n;\R^n)$,}
\]
we find by density that also for $u \in H^{s,p}(\O)$ we have
\[
\langle\mathbbm{1}_\Omega \nabla u,\psi\rangle = -\int_{\O}u \Div \psi \,dx +\int_{\partial \O} \mathrm{Tr}\,u\, \psi \cdot \nu_\O \,d\Hcal^{n-1} \quad \text{for all $\psi \in C^\infty(\R^n;\R^n)$.}
\]
Hence, by definition of the Riesz potential via duality, we find for all $\phi \in C_c^\infty(\O;\R^n)$
\[
\int_{\Omega} I_{1-s}(\mathbbm{1}_\O \nabla u) \cdot \varphi \,dx = \langle \mathbbm{1}_\O \nabla u, I_{1-s}\varphi \rangle = -\int_{\O}u \Div^s \varphi \,dx +\int_{\partial \O} \mathrm{Tr}\,u\, I_{1-s}\varphi \cdot \nu_\O \,d\Hcal^{n-1}.
\]
This proves that $G^s_\O u = I_{1-s}(\mathbbm{1}_\O \nabla u)$ and thus, $u \in \widetilde{H}^{s,p}(\O)$.

For the reverse inclusion, note that when $p>q_s'$ then $I_{1-s}^\O$ is an isomorphism from $H^{s,p'}_0(\O)$ to $W^{1,p'}(\Omega)$ by Lemma~\ref{le:Rieszisomorphism}. By density of $C_c^{\infty}(\O;\R^n)$ in $H^{s,p'}_0(\Omega;\R^n)$, we obtain
\[
\int_{\O}G^s_\O u \cdot \varphi \,dx = -\int_{\O}u\Div^s \varphi \,dx + \int_{\partial\O}gI_{1-s}\varphi \cdot \nu_\O \,d\Hcal^{n-1} \quad \text{for all $\varphi \in H^{s,p'}_0(\O;\R^n),$}
\]
and by setting $\psi = I_{1-s}^\O \varphi \in W^{1,p'}(\O;\R^n)$, we infer that
\begin{equation}\label{eq:weirdintbyparts}
\int_{\O}G^s_\O u \cdot (I_{1-s}^\O)^{-1}\psi \,dx = -\int_{\O}u\Div \psi \,dx + \int_{\partial\O}g \psi \cdot \nu_\O \,d\Hcal^{n-1} \quad \text{for all $\psi \in W^{1,p'}(\O;\R^n).$}
\end{equation}
By Lemma~\ref{le:Rieszisomorphism} and duality, we find that $I_{1-s}^\O:H^{s-1,p}_0(\O) \to L^p(\O)$ is an isomorphism, and thus, we have that $U:=(I_{1-s}^\O)^{-1}G^s_\O u \in H^{s-1,p}_0(\O;\R^n)$. By \eqref{eq:weirdintbyparts} applied to $\psi \in C_c^{\infty}(\O;\R^n)$, we find
\[
\int_{\O} U \cdot \psi \,dx = \int_{\O} G^s_\O u \cdot (I_{1-s}^\O)^{-1}\psi \,dx = -\int_{\O}u\Div \psi \,dx,
\]
which proves that $U|_{\O}=\nabla u$. Hence, $\nabla u \in H^{s-1,p}(\O;\R^n)$, which proves $u \in H^{s,p}(\O)$ by \eqref{eq:hspnorm}.

For the case $1/s < p \leq q_s'$, we may find by \eqref{eq:trace} a function $v \in H^{s,p}(\O)$ with $\Tr v = \Trt u$ given that $W^{s-1/p,p}(\partial \O)$ is the trace space of $H^{s,p}(\O)$. Combining \eqref{eq:hardytilde} with Corollary~\ref{cor:hardy}, we infer that $(u-v)/d^s \in L^p(\O)$. Moreover, since $G^s_\O v = I_{1-s}(\mathbbm{1}_\O \nabla v)$, we find
\[
\int_{\O}G^s_\O (u-v) \cdot \varphi \,dx = -\int_{\O}(u-v)\Div^s \varphi\,dx \quad \text{for all $\varphi \in C_c^{\infty}(\O;\R^n)$.}
\]
We conclude that $u-v \in \widetilde{H}^{s,p}_0(\O)$, so that the result follows by Part~$(i)$ of this theorem.
\end{proof}
\begin{remark}\label{rem:oneinclusion}
    a) The proof shows that for all $s \in (0,1)$ and $p \in (1,\infty)$
    \[
    H^{s,p}_0(\O) \subset \widetilde{H}^{s,p}_0(\O)
    \]
    and for all $sp>1$, we have that
    \[
    H^{s,p}(\O) \subset \widetilde{H}^{s,p}(\O) \quad \text{and} \quad \mathrm{Tr}=\widetilde{\mathrm{Tr}} \quad \text{on $H^{s,p}(\O)$.}
    \]

    b) In the case that \eqref{eq:weirdintbyparts} holds, then one can also conclude $u \in H^{s,p}(\O)$ for all $sp>1$, since the rest of the argument only requires $p>q_s$.

    c) The exceptional value $q_s$ arises from the fact that $H^{\frac{s-1}{2}(s)}_{q_s}(\overline{\Omega})$ is the special borderline case where both
    \[
    H^{\frac{s-1}{2}(s)}_{q_s}(\overline{\Omega}) \not = H^{s,q_s}_0(\O) \quad \text{and} \quad H^{\frac{s-1}{2}(s)}_{q_s}(\overline{\Omega})  \not \subset d^{\frac{s-1}{2}} H^{\frac{s+1}{2},q_s}(\Omega) + H^{s,q_s}_0(\Omega). \qedhere
    \]
\end{remark}

Even though the spaces only coincide for certain ranges of $s$ and $p$, it surprisingly turns out that the norms induced by the extended fractional gradients are equivalent to the norms on the Bessel potential space.

\begin{proposition}[Equivalent norms]\label{prop:equivnorm}
    Let $s \in (0,1)$, $p\in (1,\infty)$ with $p \not = q_s$, and let $\Omega \subset \R^n$ be a smooth and bounded domain. Then, we have equivalence of the norms
    \begin{equation}\label{eq:equiv0}
    \norm{u}_{L^p(\Omega)} +\norm{\bar{G}^s_\O u}_{L^p(\Omega;\R^n)} \sim \norm{u}_{H^{s,p}(\R^n)} \quad \text{for all $u \in H^{s,p}_0(\Omega)$}
    \end{equation}
    and, for $sp>1$,
    \begin{equation}\label{eq:equiv}
    \norm{u}_{L^p(\Omega)} +\norm{G^s_\O u}_{L^p(\Omega;\R^n)} \sim \norm{u}_{H^{s,p}(\Omega)} \quad \text{for all $u \in H^{s,p}(\Omega)$.}
    \end{equation}
\end{proposition}
\begin{proof}
    We first prove \eqref{eq:equiv0}. If $p<q_s$, it follows from Theorem~\ref{th:tildechar}\,$(i)$ that both norms turn $H^{s,p}_0(\O)$ into a Banach space, and since the right-hand side norm is stronger by Proposition~\ref{DsProps}, we conclude the result by the Banach isomorphism theorem. Next, let $p>q_s$ and consider $u \in C_c^{\infty}(\Omega)$. Then, we find that
    \[
    \int_{\Omega}\bar{G}^s_\O u \phi \,dx = \int_{\Omega} \nabla u I^{\Omega}_{1-s}\phi \,dx \quad \text{for all $\phi \in C_c^{\infty}(\Omega)$.}
    \]
    By density, this holds for all $\phi \in L^{p'}_0(\Omega)$, and we find for any $\psi \in H^{1-s,p'}(\O)$ that
    \begin{align*}
    \int_{\Omega} \nabla u \psi \,dx &= \int_{\Omega} \bar{G}^s_\O u(I^\Omega_{1-s})^{-1}\psi \,dx \\
    &\leq \norm{\bar{G}^s_\O u}_{L^p(\Omega;\R^n)}\norm{(I^\Omega_{1-s})^{-1}\psi}_{L^{p'}(\Omega)} \\
    &\leq C\norm{\bar{G}^s_\O u}_{L^p(\Omega;\R^n)}\norm{\psi}_{H^{1-s,p'}(\Omega)},
    \end{align*}
    where we used Lemma~\ref{le:Rieszisomorphism} in the last inequality. By \eqref{eq:hspnorm} and duality, we conclude that
    \begin{align*}
    \norm{u}_{H^{s,p}(\R^n)} &\leq C \left(\norm{u}_{L^p(\Omega)}+\norm{\nabla u}_{H^{s-1,p}(\Omega)}\right) \\
    &\leq C\left(\norm{u}_{L^p(\Omega)}+\norm{\bar{G}^s_\O u}_{L^{p}(\Omega;\R^n)}\right)
    \end{align*}
    for all $u \in C_c^{\infty}(\Omega)$. By density, this inequality holds also for all $u \in H^{s,p}_0(\Omega)$. 

    For \eqref{eq:equiv}, we note that the norm on the left is bounded by the $H^{s,p}(\O)$-norm since $G^s_\O = I_{1-s}(\mathbbm{1}_\O \nabla \cdot)|_{\O}$ on $H^{s,p}(\O)$ and \eqref{eq:estGso}. For the other inequality, let $u \in C^\infty(\overline{\Omega})$, then we compute that
    \[
    \int_{\Omega} G^s_\O u \phi\,dx = \int_{\Omega}\nabla u I^{\Omega}_{1-s}\phi \,dx \quad \text{for all $\phi \in C_c^{\infty}(\Omega)$.}
    \]
    By density, the same holds for all $\phi \in L^{p'}_0(\Omega)$, and we find for $\psi \in H^{1-s,p'}(\Omega)$ that
    \begin{align}
    \begin{split}\label{eq:estdualityGso}
    \int_{\Omega} \nabla u \psi \,dx &= \int_{\Omega} G^s_\O u(I^\Omega_{1-s})^{-1}\psi \,dx \\
    &\leq \norm{G^s_\O u}_{L^p(\Omega;\R^n)}\norm{(I^\Omega_{1-s})^{-1}\psi}_{L^{p'}(\Omega)} \\
    &\leq C\norm{G^s_\O u}_{L^p(\Omega;\R^n)}\norm{\psi}_{H^{1-s,p'}(\Omega)},
    \end{split}
    \end{align}
    where we used Lemma~\ref{le:Rieszisomorphism} in the last inequality; note that $sp>1$ implies $p>q_s$. We conclude via \eqref{eq:hspnorm} and duality that
    \begin{align*}
    \norm{u}_{H^{s,p}(\Omega)} &\leq C \left(\norm{u}_{L^p(\Omega)}+\norm{\nabla u}_{H^{s-1,p}(\Omega)}\right) \\
    &\leq C\left(\norm{u}_{L^p(\Omega)}+\norm{G^s_\O u}_{L^{p}(\Omega;\R^n)}\right)
    \end{align*}
    for all $u \in C^{\infty}(\overline{\Omega})$. By density, this inequality holds also for all $u \in H^{s,p}(\Omega)$, which finishes the proof.
\end{proof}
This equivalence enables us to also prove fractional Poincar\'e-Wirtinger inequalities. In fact, for certain values of $s$ and $p$ we need not even subtract mean values, due to the following observation.
\begin{lemma}\label{le:zerofunctions}
    Let $s \in (0,1)$, $p\in (1,\infty)$ with $p>q_s$ and $\Omega \subset \R^n$ be a smooth and bounded domain. Then, if $u \in H^{s,p}_0(\O)$ satisfies $\bar{G}^s_\O u = 0$ on $\O$, we have that $u = 0$. If instead $p \leq q_s$, the same holds when $(u)_\O=0$.
\end{lemma}
\begin{proof}
    Any such function will be given by $u=(I_{1-s}^{\O})^{-1} c$ for some $c \in \R$. This constitutes a one-dimensional linear space, and we will show that any non-zero element, i.e., $c \not =0$, does not lie in $H^{s,p}_0(\O)$ when $p>q_s$. This follows by showing that $u$ has blow-up behavior near the boundary, and thus, does not satisfy the Hardy inequality.  

    To prove this, consider the energy
    \[
    \Ical(\mu):= \int_{\R^n}\int_{\R^n} |y-x|^{-(n+s-1)}\,d\mu(y)\,d\mu(x) \quad \text{for $\mu \in \Pcal(\overline{\Omega})$}
    \]
    where $\Pcal(\overline{\Omega})$ denotes the probability measures supported on $\overline{\Omega}$. It follows from \cite[Theorem~4.1 and 4.2]{HRS19} that $\Ical$ admits a unique minimizer $\mu_0$, which satisfies
    \[
    I_{1-s}\mu_0 (x) = c \quad \text{on $\Omega$,}
    \]
    for some non-zero constant $c>0$. Note also that
    \[
    \int_{\R^n}\int_{\R^n} |y-x|^{-(n+s-1)}\,d\mu_0(y)\,d\mu_0(x) \propto \int_{\R^n}  |\xi|^{s-1}|\widehat{\mu_0}(\xi)|^2\,d\xi,
    \]
    with $\widehat{\mu_0} \in C^\infty(\R^n;\C)$ the Fourier transform of $\mu_0$; this is simply \eqref{eq:rieszidentity} extended to probability measures. In particular, $\mu_0 \in H^{\frac{s-1}{2}}_0(\Omega)$, so that we may apply Lemma~\ref{le:Rieszisomorphism0} to conclude that $\mu_0$ must agree with $u=(I_{1-s}^{\O})^{-1} c$. In addition, since the constant function is smooth, we can use Theorem~\ref{th:regularity} to infer that $u \in \bigcap_{m \in \R} H^{\frac{s-1}{2}(m)}_2(\overline{\Omega})$, that is, $u=d^{(s-1)/2}v$ for some $v \in C^{\infty}(\overline{\Omega})$, cf.~\cite[Eq.(4.1)]{Gru15}. If we now consider the function $V:=I_{1-s}u$, then it is continuous because it is the convolution of a locally integrable function with a compactly supported integrable function. Moreover, we can compute that
    \[
    V = c \quad \text{in $\overline{\Omega}$} \quad \text{and} \quad \Delta V = -(-\Delta)^{\frac{s+1}{2}}u > 0 \quad \text{in $\overline{\Omega}^c$,}
    \]
    where the second part uses that $u$ is non-negative and not identically zero, since it agrees with the probability measure $\mu_0$. If we now use the fact that $V$ vanishes at infinity and is smooth outside $\overline{\Omega}$, we infer from the maximum principle that $V$ attains its maximum exclusively on $\overline{\O}$. The Hopf lemma then implies that the derivative of $V$ on $\partial \Omega$ in the direction of the inner unit normal is strictly positive, at points where it exists \cite[Lemma~3.4]{GiT01}.

    On the other hand, if $v(x_0)=0$ for some point $x_0 \in \partial \Omega$, it can be shown that $V$ is differentiable there, and thus, its derivative is zero given that $V$ is constant on $\Omega$. Indeed, we then have the estimate $|u(y)|\leq C|x_0-y|d(y)^{(s-1)/2}$, so its convolution with the differentiated kernel of the Riesz potential converges at $x_0$. Since this would contradict the strict positivity of the normal derivative, we conclude that $v$ is strictly positive on $\partial \O$. In particular, we deduce that $u/d^s$ behaves like $d^{-(1+s)/2}$ near the boundary. Therefore, $u/d^s \not \in L^p(\Omega)$ for $p>q_s$, which shows that $u \not \in H^{s,p}_0(\Omega)$ by the fractional Hardy inequality \cite[Proposition~I.5.7]{Tri01}.

    The case when $p\leq q_s$ and the mean vanishes $(u)_\Omega =0$ follows immediately from the non-negativity of $u$.
\end{proof}
\begin{remark}\label{rem:singularities}
    We note that the functions for which $\bar{G}^s_\O u=0$ lie in $\widetilde{H}^{s,p}_0(\O)$ for $p<q_s'$, since they can be written as $u=d^{(s-1)/2}v$ for some $v \in C^{\infty}(\overline{\Omega})$, which is contained in $L^p(\O)$. However, the previous result shows that they only lie in $H^{s,p}_0(\O)$ when also $p<q_s$. In particular, in light of the equivalence of norms from Proposition~\ref{prop:equivnorm}, these functions cannot be approximated by $C_c^{\infty}(\O)$ functions in the natural norm $\|\cdot\|_{L^p(\O)}+\|\bar{G}^s_\O \cdot \|_{L^p(\O;\R^n)}$ when $p \in (q_s,q_s')$, showing that $\widetilde{H}^{s,p}_0(\O)$ cannot be defined via density with respect to this norm. This peculiar fact was also recently noticed in \cite{FLP26} for the case $p=2$.
\end{remark}
We obtain the following Poincar\'e-Wirtinger inequalities, which open up the possibility for using the extended fractional gradients in applications of partial differential equations and variational problems in domains.
\begin{corollary}[Poincar\'e-Wirtinger inequalities]\label{cor:poincwirt}
     Let $s \in (0,1)$, $p\in (1,\infty)$ and $\Omega \subset \R^n$ be a smooth and bounded domain. 
     \begin{itemize}
         \item[(i)] If $p<q_s$, there exists a constant $C>0$ such that
     \[
     \|u-\mathbbm{1}_\O(u)_\O\|_{H^{s,p}(\R^n)} \leq C \|\bar{G}^s_\O u\|_{L^p(\O;\R^n)} \quad \text{for all $u \in H^{s,p}_0(\O)$.}
     \]
     If $p>q_s$, there exists a constant $C>0$ such that
     \[
     \|u\|_{H^{s,p}(\R^n)} \leq C \|\bar{G}^s_\O u\|_{L^p(\O;\R^n)} \quad \text{for all $u \in H^{s,p}_0(\O)$.}
     \] \\[-0.6cm]
     
     \item[(ii)] If $sp>1$, there exists a $C>0$ such that
     \[
     \|u-(u)_\O\|_{H^{s,p}(\O)} \leq C \|G^s_\O u\|_{L^p(\O;\R^n)} \quad \text{for all $u \in H^{s,p}(\O)$.}
     \]
     \end{itemize}
\end{corollary}
\begin{proof}
    Part\,$(i)$ follows from Proposition~\ref{prop:equivnorm} and Lemma~\ref{le:zerofunctions} using the usual contradiction argument. For Part\,$(ii)$, suppose that $G^s_\O u = 0$ on $\Omega$. Then, $\nabla u =0$ on $\Omega$ by \eqref{eq:estdualityGso}. We conclude that $u$ must be constant, so we can again use the usual contradiction argument.
\end{proof}

The results from this section allow us to characterize the Bessel potential spaces via the regional fractional gradients.

\begin{proof}[Proof of Theorem~\ref{th:char}]
Note that the one implication is clear by Theorem~\ref{th:emb1} together with the fractional Hardy inequality (\cite[Proposition~I.5.7]{Tri01}) and Corollary~\ref{cor:hardy}. For the reverse, we distinguish the two cases.

    \textit{Part (i).} We note that by \eqref{eq:connDsobarGso}
    \[
    \bar{G}^s_\O u=D^s_\O u+uD^s \mathbbm{1}_\O \in L^p(\Omega;\R^n),
    \]
    since $|D^s \mathbbm{1}_\O| \leq C d^{-s}$, cf.~\eqref{eq:distest}. Hence, the result follows by Theorem~\ref{th:tildechar} and Proposition~\ref{prop:equivnorm}. \medskip

    \textit{Part (ii).} Since $W^{s-1/p,p}(\partial \O)$ is the trace space of $H^{s,p}(\O)$ by \eqref{eq:trace}, there exists a $v \in H^{s,p}(\O)$ with its trace equal to $g$. Therefore, by Corollary~\ref{cor:hardy}, we deduce that
    \[
    \int_{\partial \Omega}\nu_\Omega(y) \frac{g(y)-v(\cdot)}{|y-\cdot|^{n+s-1}}\,d\Hcal^{n-1}(y) \in L^p(\O;\R^n),
    \]
    and thus,
    \[
    (u-v)/d^s \in L^p(\O).
    \]
    Since also $v \in \Hcal^{s,p}(\O)$ by Theorem~\ref{th:emb1}, we infer $u-v \in H^{s,p}_0(\O)$ by Part\,$(i)$. Hence, $u \in H^{s,p}(\O)$ and $\mathrm{Tr}\,u=\mathrm{Tr}\,v=g$. Finally, by combining \eqref{eq:connDsoGso2} and \eqref{eq:connDsoGso3}, and extending it weakly, one finds
    \[
    G^s_\O u =D^s_\O u- \frac{1}{\gamma_{n,1-s}}\int_{\partial \O} \nu_\O(y)\frac{\mathrm{Tr}\,u(y)-u(\cdot)}{|y-\cdot|^{n+s-1}}\,d\Hcal^{n-1}(y)
    \]
    with $\gamma_{n,1-s}$ the normalizing constant of the Riesz potential. Hence, the inequality of the norms follows from Corollary~\ref{cor:poincwirt} and the triangle inequality.
\end{proof}
\begin{remark}
    a) Let $p>1/(1-s)$ and $u \in L^p(\O)$, then it follows from H\"older's inequality that
    \[
    u/d^s \in L^r(\O) \quad \text{for all $r<\frac{p}{1+ps}$.}
    \]
    In particular, it then holds that
    \[
    \Hcal^{s,p}(\O) \hookrightarrow H^{s,r}_0(\O),
    \]
    and thus, via an interpolation inequality, $\Hcal^{s,p}(\O)$ is compactly embedded into $L^q(\O)$ for $q<p$; this compact embedding actually holds for all $p \in (1,\infty)$ using Theorem~\ref{th:regional-local-bessel} and that $L^p(\O) \cap H^{s,p}_{\rm loc}(\O)$ is compactly embedded into $L^q(\O)$ for $q<p$. \medskip
    
    b) It is still an open problem to see whether
    \begin{equation}\label{eq:conjecture}
    H^{s,p}(\O)= \Hcal^{s,p}(\O).
    \end{equation}
    Indeed, any function $u \in \widetilde{H}^{s,p}_0(\O) \setminus H^{s,p}_0(\O)$ cannot satisfy the Hardy inequality by Theorem~\ref{th:tildechar}, which implies
    \[
    D^s_\O u =\bar{G}^s_\O u -uD^s \mathbbm{1}_\O \not \in L^p(\O;\R^n),
    \]
    that is, $u \not \in \Hcal^{s,p}(\O)$. Therefore, these functions do not provide a counterexample to the inclusion $\Hcal^{s,p}(\O) \subset H^{s,p}(\O)$. Moreover, by Theorem~\ref{th:char}, it would be sufficient to prove Hardy inequalities on the space $\Hcal^{s,p}(\O)$ in order to conclude \eqref{eq:conjecture}.
\end{remark}

\subsection*{Acknowledgements}
This work was supported by {\it Agencia Estatal de Investigaci\'on} (Spain) through grant PID2023-151823NB-I00 and {\it Junta de Comunidades de Castilla-La Mancha} (Spain) through grant SBPLY/23/180225/000023. G.G.-S. is supported by a Doctoral Fellowship by Universidad de Castilla-La Mancha \text{2024-UNIVERS-12844-404}. H.S. acknowledges support by the NWO through the Vidi grant \textit{SPARGO: Exploring and Exploiting the Geometric Landscape of Infinite-Dimensional Sparse Optimization} (Grant Number VI.Vidi.243.200) and by the Fonds de la Recherche Scientifique - FNRS through a MIS-Ulysse project (Scientific Impulse Mandate Instrument) number F.6002.25.

\subsection*{Generative AI declaration}
OpenAI's ChatGPT was used as an interactive mathematical assistant during the preparation of this article. Beyond general expository improvements, it suggested to us the main idea of Proposition~\ref{prop:nullspace} consisting of composing with a convex function and testing with the distance function to the boundary. It also pointed us to the reference~\cite{li2019} for obtaining the proof of Theorem~\ref{th:emb1}, and constructed the main idea behind~Lemma~\ref{le:zerofunctions}, in particular, the connection with equilibrium measures in \cite{HRS19} and the argument with Hopf's lemma. All other mathematical ideas and results were developed by the authors; every AI-assisted argument and calculation was independently reconstructed and verified by the authors, who take full responsibility for the final manuscript.

\addcontentsline{toc}{section}{\protect\numberline{}References}
\bibliographystyle{abbrv}
\bibliography{bibliography}
\end{document}